# Robust Sum MSE Optimization for Downlink Multiuser MIMO Systems with Arbitrary Power Constraint: Generalized Duality Approach

Tadilo Endeshaw Bogale, *Student Member, IEEE* and Luc Vandendorpe *Fellow, IEEE*

*Abstract*— This paper considers linear minimum mean-square-error (MMSE) transceiver design problems for downlink multiuser multiple-input multiple-output (MIMO) systems where imperfect channel state information is available at the base station (BS) and mobile stations (MSs). We examine robust sum mean-square-error (MSE) minimization problems. The problems are examined for the generalized scenario where the power constraint is per BS, per BS antenna, per user or per symbol, and the noise vector of each MS is a zero-mean circularly symmetric complex Gaussian random variable with arbitrary covariance matrix. For each of these problems, we propose a novel duality based iterative solution. Each of these problems is solved as follows. First, we establish a novel sum average mean-square-error (AMSE) duality. Second, we formulate the power allocation part of the problem in the downlink channel as a Geometric Program (GP). Third, using the duality result and the solution of GP, we utilize alternating optimization technique to solve the original downlink problem. To solve robust sum MSE minimization constrained with per BS antenna and per BS power problems, we have established novel downlink-uplink duality. On the other hand, to solve robust sum MSE minimization constrained with per user and per symbol power problems, we have established novel downlink-interference duality. For the total BS power constrained robust sum MSE minimization problem, the current duality is established by modifying the constraint function of the dual uplink channel problem. And, for the robust sum MSE minimization with per BS antenna and per user (symbol) power constraint problems, our duality are established by formulating the noise covariance matrices of the uplink and interference channels as fixed point functions, respectively. We also show that our sum AMSE duality are able to solve other sum MSE-based robust design problems. Computer simulations verify the robustness of the proposed robust designs compared to the non-robust/naive designs.

## I. INTRODUCTION

The spectral efficiency of wireless channels can be enhanced by utilizing multiple-input multiple-output (MIMO) systems. This performance improvement is achieved by exploiting the transmit and receive diversity. In [1], a fundamental relation between mutual information and minimum mean-square-error (MMSE) has been derived for MIMO Gaussian channels. Furthermore, it has been shown that different transceiver optimization problems are equivalently formulated as a function of MMSE matrix, for instance, minimizing bit error rate, maximizing capacity etc [2]–[4]. For these reasons, mean-square-error (MSE) based problems are commonly examined in multiuser networks.

The linear MMSE transceiver design problems for multiuser MIMO systems can be examined in the uplink and/or downlink channels. In [5] and [6], sum MSE minimization problem is examined in the uplink channel. The authors of these papers exploit the fact that the sum MSE optimization problem in the uplink channel (after rank relaxation) can be formulated as a convex optimization problem for which global optimal solution can be obtained efficiently. It is well know that transceiver design problems in the uplink channel are better understood than that of the downlink channel. Due to this fact most literatures (also our current paper) focus on solving transceiver design problems in the downlink channel. In [7], the sum MSE minimization problem is addressed in the downlink channel. This paper has solved the latter problem directly in the downlink channel. In this channel, however, since the precoders of all users are jointly coupled, downlink MSE-based transceiver design problems have more complicated mathematical structure than their dual uplink problems [3], [8].

In [3] and [9], the sum MSE minimization constrained with a total base station (BS) power problem is considered in the downlink channel. The authors of these papers solve the downlink sum MSE minimization problem by examining the equivalent uplink problem, and then establishing the MSE duality between uplink and downlink channels. These two papers show that duality based approach of solving the downlink problem has easier to handle mathematical structure than that of [7]. Moreover, in some cases, duality solution approach can exploit the hidden convexity of the downlink channel MSE-based problems (for example minimization of sum MSE constrained with a total BS power problem [3], [8]). These duality are established by assuming that perfect channel state information (CSI) is available at the BS and mobile stations (MSs). However, due to the inevitability of channel estimation error, CSI can never be perfect. This motivates [10] to establish the MSE duality under imperfect CSI for multiple-input single-output (MISO) systems. The latter work is extended in [11] for MIMO case. In [12], the MSE downlink-uplink duality has been established by considering imperfect CSI both at the BS

The authors would like to thank the Region Wallonne for the financial support of the project MIMOCOM in the framework of which this work has been achieved. Part of this work has been presented in the 45th Annual Conference on Information Sciences and Systems (CISS), Baltimore, Maryland, Mar. 2011, and in the 22nd IEEE International Symposium on Personal, Indoor and Mobile Radio Communications (PIMRC), Toronto, Canada, Sep. 2011. Tadilo E. Bogale and Luc Vandendorpe are with the ICTEAM Institute, Université catholique de Louvain, Place du Levant 2, 1348 - Louvain La Neuve, Belgium. Email: {tadilo.bogale, luc.vandendorpe}@uclouvain.be, Phone: +3210478071, Fax: +3210472089.

and MSs, and with antenna correlation only at the BS. This duality is examined by analyzing the Karush-Kuhn-Tucker conditions of the uplink and downlink channel problems. In [13], we have established three kinds of MSE uplink-downlink duality by considering that imperfect CSI is available both at the BS and MSs, and with antenna correlation only at the BS. These duality are established by extending the three level MSE duality of [8] to imperfect CSI. In [14], we have shown that the three kinds of MSE duality known from [8] can be established for the perfect and imperfect CSI scenario just by transforming the power allocation matrices from uplink to downlink channel and vice versa. All of these MSE duality are established by assuming that the entries of the noise vector of each MS are independent and identically distributed (i.i.d) zero-mean circularly symmetric complex Gaussian (ZMCSCG) random variables all with the same variance. However, in practice, MSs are spaced far apart from each other and the noise vector of each MS may include other interference Gaussian signals [2]. For these reasons, the noise variances of all MSs are not necessarily the same. This motivates [15] to exploit the three kinds of MSE downlink-uplink duality known in [8] and [14] for the scenario where the noise vector of each MS is a ZMCSCG random variable with arbitrary covariance matrix. However, still the work of [15] exploits the MSE downlink-uplink duality by assuming that perfect CSI is available at the BS and MSs. Moreover, the duality of all of the aforementioned papers including [15] are able to solve total BS power constrained MSE or average mean-square-error (AMSE)-based problems only. In [16] (see also [17] and [18]), duality based iterative solutions for rate and signal-to-interference-and-noise ratio (SINR)-based problems with multiple linear transmit covariance constraints (LTCC) have been proposed. The problems of [16] are examined as follows. First, new auxiliary variables are incorporated to combine the LTCC into one constraint. Second, keeping the introduced auxiliary variables constant, the existing duality approach is applied to get the optimal covariance matrices. Third, keeping the covariance matrices of all users constant, the introduced variables are updated by sub-gradient optimization method. Fourth, the second and third steps are repeated until convergence. However, the aforementioned iterative algorithm has one major drawback. For fixed introduced variables, if the global optimality of the downlink (dual uplink) problem is not ensured, this iterative algorithm is not guaranteed to converge[1]. Besides this drawback, the application of the latter duality for MSE-based problems is not clear.

In a practical multi-antenna BS system, the maximum power of each BS antenna is limited [19]. Due to this, we developed duality based iterative algorithms to solve sum MSE-based constrained with each BS antenna power design problems for the perfect CSI scenario in [20]. On the other hand, channel estimation error is inevitable, and in some scenario allocating different powers to different users or symbols according to their channel conditions (to ensure fairness among users or symbols) has practical interest. This motivates us to first propose novel generalized sum AMSE duality and then we utilize our new duality results to examine the following sum MSE-based robust design problems. Robust sum MSE minimization with a per total BS ($\mathcal{P}1$), per BS antenna ($\mathcal{P}2$), per user ($\mathcal{P}3$) and per symbol ($\mathcal{P}4$) power constraints, respectively. These problems are examined by considering that the BS and MS antennas exhibit spatial correlations, and the CSI at both ends is imperfect. The robustness against imperfect CSI is incorporated into our designs using stochastic approach [14]. Moreover, the noise vector of each MS is a ZMCSCG random variable with arbitrary covariance matrix.

To the best of our knowledge, the problems $\mathcal{P}1 - \mathcal{P}4$ are not convex. Furthermore, duality based solutions for these problems with our noise covariance matrix assumptions are not known. In the current paper, we propose duality based iterative solutions to solve the problems. Each of these problems is solved as follows. First, we establish novel sum AMSE duality. Second, we formulate the power allocation part of the problem in the downlink channel as a Geometric Program (GP). Third, using the duality result and the solution of GP, we utilize alternating optimization technique to solve the original downlink problem. We have established novel downlink-uplink duality to solve $\mathcal{P}1$ and $\mathcal{P}2$. On the other hand, we have established novel downlink-interference duality to solve $\mathcal{P}3$ and $\mathcal{P}4$. For problem $\mathcal{P}1$, the current duality is established by modifying the constraint function of the dual uplink channel problem. And, for the problems $\mathcal{P}2$ and $\mathcal{P}3 - \mathcal{P}4$, our duality are established by formulating the noise covariance matrices of the uplink and interference channels as fixed point functions, respectively. The duality of this paper generalize all existing sum MSE (AMSE) duality. The main contributions of the current paper can thus be summarized as follows.

1) We have established novel sum AMSE downlink-uplink and downlink-interference duality to solve the problems $\mathcal{P}1 - \mathcal{P}4$ for the generalized scenario where the noise vector of each MS is a ZMCSCG random variable with arbitrary covariance matrix, the BS and MS antennas exhibit spatial correlations and the CSI at both ends is imperfect. The downlink-uplink duality are established to solve $\mathcal{P}1$ and $\mathcal{P}2$, whereas the downlink-interference duality are established to solve $\mathcal{P}3$ and $\mathcal{P}4$. For problem $\mathcal{P}1$, the current duality is established by modifying the constraint function of the dual uplink channel problem. And, for the problems $\mathcal{P}2$ and $\mathcal{P}3 - \mathcal{P}4$, our duality are established by formulating the noise covariance matrices of the uplink and interference channels as fixed point functions, respectively. Therefore, the current sum AMSE duality generalize the hitherto sum AMSE duality.

2) By employing the system model of [4] and [15], we formulate the power allocation part of our robust sum MSE minimization problems in the downlink channel as GPs. These GPs are formulated by modifying the GP formulation approach of [13] to our scenario. Consequently, our GP solution and duality results enable us to solve the problems $\mathcal{P}1$ - $\mathcal{P}4$ by applying alternating optimization techniques (i.e., duality based iterative algorithms) of [3] and [14]. As will be clear later, the current duality based iterative algorithms can solve other sum MSE-based robust design problems.
3) In our simulation results, we have observed that the proposed duality based iterative algorithms utilize less total BS power than that of existing algorithms for the problems $\mathcal{P}2 - \mathcal{P}4$.
4) We examine the effects of channel estimation errors and antenna correlations on the system performance.

---

[1] Note that for fixed introduced variables, the global optimum covariance matrices for the problems of [16] can be obtained. Thus, the duality-based iterative algorithm of this paper is guaranteed to converge for the problems of [16]. However, as will be clear later, since all of our problems are non-convex, the iterative algorithm of [16] is not guaranteed to converge for our problems.

The remaining part of this paper is organized as follows. In Section II, multiuser MIMO downlink, and virtual uplink and interference system models are presented. In Section III, MIMO channel model under imperfect CSI is discussed. In Section IV, we formulate our robust design problems $\mathcal{P}1$ - $\mathcal{P}4$ and discuss the general framework of our duality based iterative solutions. Sections V - IX present the proposed duality based iterative solutions for solving the problems $\mathcal{P}1$ - $\mathcal{P}4$. The extension of our duality based iterative algorithms to other problems is discussed in Section X. In Section XI, computer simulations are used to compare the performance of the proposed duality based iterative algorithms with that of existing algorithms, and the robust and non-robust/naive designs. Finally, conclusions are drawn in Section XII.

*Notations:* The following notations are used throughout this paper. Upper/lower case boldface letters denote matrices/column vectors. The $[\mathbf{X}]_{(n,n)}$, $[\mathbf{X}]_{(n,:)}$, $\text{tr}(\mathbf{X})$, $\mathbf{X}^T$, $\mathbf{X}^H$ and $\mathrm{E}(\mathbf{X})$ denote the $(n, n)$ element, $n$th row, trace, transpose, conjugate transpose and expected value of $\mathbf{X}$, respectively. $\mathbf{I}_n(\mathbf{I})$ is an identity matrix of size $n \times n$ (appropriate size) and $\mathcal{C}(\Re)^{M \times M}$ represent spaces of $M \times M$ matrices with complex (real) entries. The diagonal and block-diagonal matrices are represented by $\text{diag}(.)$ and $\text{blkdiag}(.)$ respectively. Subject to is denoted by s.t and $(.)^\star$ denotes optimal solution. The $n$th norm of a vector $\mathbf{x}$ is represented by $\|\mathbf{x}\|_n$. The superscripts $(.)^{DL}$, $(.)^{UL}$ and $(.)^I$ denotes downlink, uplink and interference, respectively.

## II. SYSTEM MODEL

In this section, multiuser MIMO downlink, and virtual uplink and interference system models are considered. The BS equipped with $N$ transmit antennas is serving $K$ MSs. Each MS has $M_k$ antennas to multiplex $S_k$ symbols. The total number of MS antennas and symbols are $M = \sum_{k=1}^{K} M_k$ and $S = \sum_{k=1}^{K} S_k$, respectively. The entire symbol can be written in a data vector $\mathbf{d} = [\mathbf{d}_1^T, \cdots, \mathbf{d}_K^T]^T$, where $\mathbf{d}_k \in \mathcal{C}^{S_k \times 1}$ is the symbol vector for the $k$th MS. In the downlink channel (see Fig. 1.(a)), the BS precodes $\mathbf{d} \in \mathcal{C}^{S \times 1}$ into an $N$ length vector by using the overall precoder matrix $\mathbf{B} = [\mathbf{B}_1, \cdots, \mathbf{B}_K]$, where $\mathbf{B}_k \in \mathcal{C}^{N \times S_k}$ is the precoder matrix for the $k$th MS. The $k$th MS uses a linear receiver $\mathbf{W}_k \in \mathcal{C}^{M_k \times S_k}$ to recover its symbol $\mathbf{d}_k$ as

$$\widehat{\mathbf{d}}_k^{DL} = \mathbf{W}_k^H (\mathbf{H}_k^H \mathbf{B} \mathbf{d} + \mathbf{n}_k) = \mathbf{W}_k^H (\mathbf{H}_k^H \sum_{i=1}^{K} \mathbf{B}_i \mathbf{d}_i + \mathbf{n}_k)$$

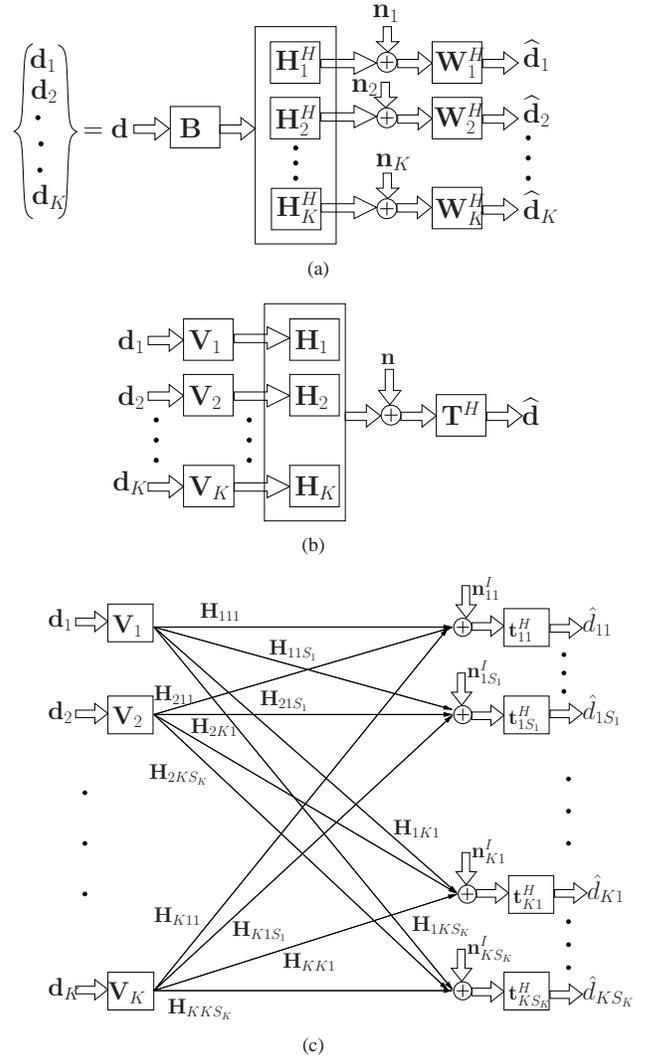

Fig. 1. Multiuser MIMO system model. (a) downlink channel. (b) virtual uplink channel. (c) virtual interference channel.

where $\mathbf{n}_k \in \mathcal{C}^{M_k \times 1}$ is the additive Gaussian noise at the $k$th MS and $\mathbf{H}_k^H \in \mathcal{C}^{M_k \times N}$ is the MIMO channel between the BS and $k$th MS. Without loss of generality, we can assume that the entries of $\mathbf{d}_k$ are i.i.d ZMCSCG random variables all with unit variance[2], i.e., $\mathrm{E}\{\mathbf{d}_k \mathbf{d}_k^H\} = \mathbf{I}_{S_k}$, $\mathrm{E}\{\mathbf{d}_k \mathbf{d}_i^H\} = \mathbf{0}$, $\forall i \neq k$, and $\mathrm{E}\{\mathbf{d}_k \mathbf{n}_i^H\} = \mathbf{0}$. The noise vector of the $k$th MS is a ZMCSCG random variable with positive semi-definite covariance matrix $\mathbf{R}_{nk} \in \mathcal{C}^{M_k \times M_k}$. To establish the sum AMSE duality, we model the virtual uplink and interference channels as shown in Fig. 1.(b)-(c). The virtual uplink and interference channels are modeled by introducing precoders $\{\mathbf{V}_k = [\mathbf{v}_{k1}, \cdots, \mathbf{v}_{kS_k}]\}_{k=1}^{K}$ and decoders $\{\mathbf{T}_k = [\mathbf{t}_{k1}, \cdots, \mathbf{t}_{kS_k}]\}_{k=1}^{K}$, where $\mathbf{v}_{ks} \in \mathcal{C}^{M_k \times 1}$ and $\mathbf{t}_{ks} \in \mathcal{C}^{N \times 1}, \forall k$. In the uplink channel, it is assumed that $\mathbf{n}$ is a ZMCSCG random variable with diagonal covariance matrix $\boldsymbol{\Psi} \in \Re^{N \times N} = \text{diag}(\psi_1, \cdots, \psi_N)$. In the interference channel, it is assumed that the $k$th user's $s$th symbol is estimated independently by $\mathbf{t}_{ks} \in \mathcal{C}^{N \times 1}$. Moreover, $\{\mathbf{n}_{ks}^I, \forall s\}_{k=1}^K$ (Fig. 1.c) are also ZMCSCG random variables with covariance matrices $\{\boldsymbol{\Delta}_{ks} \in \Re^{N \times N} = \text{diag}(\delta_{ks1}, \cdots, \delta_{ksN}), \forall s\}_{k=1}^K$ and the channels between the $k$th transmitter and all receivers are the same (i.e., $\{\mathbf{H}_{kjs} = \mathbf{H}_k, \forall j, s\}_{k=1}^K$).

---
[2]We would like to mention here that when the symbols are not white pre-whitening operation can be applied before the transmitter (precoder) and the inverse operation can be performed after the receiver (decoder) [2].

In both the uplink and interference system models, the transmitters are the same (i.e., $\{\mathbf{V}_k\}_{k=1}^K$). Furthermore, the channels between all transmitters and each receiver (one receiver in the case of uplink system model and $S$ receivers in the case of interference system model) are the same (i.e., $\{\mathbf{H}_k\}_{k=1}^K$). As will be clear later, the virtual uplink and interference system models differ on their noise covariance matrices (which is clearly seen from Fig. 1.(b)-(c)) and the approach where the decoders are designed. In the case of downlink-uplink duality, all decoders of the uplink channel are designed with a common noise covariance matrix $\boldsymbol{\Psi}$, whereas, in the case of downlink-interference duality, the $k$th user $s$th symbol decoder of the interference channel $\mathbf{t}_{ks}$ is designed by employing its own noise covariance matrix $\boldsymbol{\Delta}_{ks}$.

## III. CHANNEL MODEL

Considering antenna correlation at the BS and MSs, we model the Rayleigh fading MIMO channels as $\mathbf{H}_k^H = \widetilde{\mathbf{R}}_{mk}^{1/2} \mathbf{H}_{wk}^H \mathbf{R}_{bk}^{1/2}, \forall k$, where the elements of $\{\mathbf{H}_{wk}^H\}_{k=1}^K$ are i.i.d ZMCSCG random variables all with unit variance, and $\mathbf{R}_{bk} \in C^{N \times N}$ and $\widetilde{\mathbf{R}}_{mk} \in C^{M_k \times M_k}$ are antenna correlation matrices at the BS and MSs, respectively [14], [21]. The channel estimation is performed on $\{\mathbf{H}_{wk}^H\}_{k=1}^K$ using an orthogonal training method [22]. Upon doing so, the $k$th user true channel $\mathbf{H}_k^H$ and its minimum-mean-square-error (MMSE) estimate $\widehat{\mathbf{H}}_k^H$ can be related as (see Section II.B of [22] for more details about the channel estimation process)

$$\mathbf{H}_k^H = \widehat{\mathbf{H}}_k^H + \mathbf{R}_{mk}^{1/2} \mathbf{E}_{wk}^H \mathbf{R}_{bk}^{1/2} = \widehat{\mathbf{H}}_k^H + \mathbf{E}_k^H, \; \forall k \quad (1)$$

where $\mathbf{R}_{mk} = (\mathbf{I}_{M_k} + \sigma_{ek}^2 \widetilde{\mathbf{R}}_{mk}^{-1})^{-1}$, $\mathbf{E}_k^H$ is the estimation error and the entries of $\mathbf{E}_{wk}^H$ are i.i.d with $\mathcal{CN}(0, \sigma_{ek}^2)$. By robustness, we mean that $\{\mathbf{E}_{wk}^H\}_{k=1}^K$ are unknown but $\{\widehat{\mathbf{H}}_k^H, \mathbf{R}_{bk}, \widetilde{\mathbf{R}}_{mk}$ and $\sigma_{ek}^2\}_{k=1}^K$ are known. We assume that each MS estimates its channel and feeds the estimated channel back to the BS without any error and delay. Thus, both the BS and MSs have the same channel imperfections. With these assumptions, the downlink instantaneous MSE ($\boldsymbol{\xi}_k^{DL}$) and AMSE ($\bar{\boldsymbol{\xi}}_k^{DL}$) matrices of the $k$th user are given by

$$\boldsymbol{\xi}_k^{DL} = \mathbb{E}_{\mathbf{d}, \mathbf{n}_k} \{(\mathbf{d}_k - \widehat{\mathbf{d}}_k^{DL})(\mathbf{d}_k - \widehat{\mathbf{d}}_k^{DL})^H\}$$
$$= \mathbf{I}_{S_k} + \mathbf{W}_k^H (\mathbf{H}_k^H \sum_{i=1}^K \mathbf{B}_i \mathbf{B}_i^H \mathbf{H}_k + \mathbf{R}_{nk}) \mathbf{W}_k -$$
$$\mathbf{W}_k^H \mathbf{H}_k^H \mathbf{B}_k - \mathbf{B}_k^H \mathbf{H}_k \mathbf{W}_k$$
$$\bar{\boldsymbol{\xi}}_k^{DL} = \mathbb{E}_{\mathbf{E}_{wk}^H} \{\boldsymbol{\xi}_k^{DL}\} = \mathbf{I}_{S_k} + \mathbf{W}_k^H \boldsymbol{\Gamma}_k^{DL} \mathbf{W}_k - \mathbf{W}_k^H \widehat{\mathbf{H}}_k^H \mathbf{B}_k -$$
$$\mathbf{B}_k^H \widehat{\mathbf{H}}_k \mathbf{W}_k \quad (2)$$

where $\boldsymbol{\Gamma}_k^{DL} = \widehat{\mathbf{H}}_k^H \mathbf{B} \mathbf{B}^H \widehat{\mathbf{H}}_k + \sigma_{ek}^2 \text{tr}\{\mathbf{R}_{bk} \mathbf{B} \mathbf{B}^H\} \mathbf{R}_{mk} + \mathbf{R}_{nk}$. The total sum AMSE $\bar{\xi}^{DL} = \sum_{k=1}^K \text{tr}\{\bar{\boldsymbol{\xi}}_k^{DL}\}$ is given by

$$\bar{\xi}^{DL} = S + \sum_{k=1}^K \text{tr}\{\mathbf{W}_k^H \boldsymbol{\Gamma}_k^{DL} \mathbf{W}_k - \mathbf{W}_k^H \widehat{\mathbf{H}}_k^H \mathbf{B}_k -$$
$$\mathbf{B}_k^H \widehat{\mathbf{H}}_k \mathbf{W}_k\}. \quad (3)$$

## IV. PROBLEM FORMULATIONS

Mathematically, the aforementioned robust sum MSE minimization problems can be formulated as

$$\mathcal{P}1 : \min_{\{\mathbf{B}_k, \mathbf{W}_k\}_{k=1}^K} \bar{\xi}^{DL}, \; \text{s.t} \; \text{tr}(\sum_{k=1}^K \mathbf{B}_k \mathbf{B}_k^H) \le P_{max} \quad (4)$$

$$\mathcal{P}2 : \min_{\{\mathbf{B}_k, \mathbf{W}_k\}_{k=1}^K} \bar{\xi}^{DL}, \; \text{s.t} \; [\sum_{k=1}^K \mathbf{B}_k \mathbf{B}_k^H]_{n,n} \le \check{p}_n, \forall n \quad (5)$$

$$\mathcal{P}3 : \min_{\{\mathbf{B}_k, \mathbf{W}_k\}_{k=1}^K} \bar{\xi}^{DL}, \; \text{s.t} \; \text{tr}\{\mathbf{B}_k \mathbf{B}_k^H\} \le \check{p}_k, \forall k \quad (6)$$

$$\mathcal{P}4 : \min_{\{\mathbf{B}_k, \mathbf{W}_k\}_{k=1}^K} \bar{\xi}^{DL}, \; \text{s.t} \; \mathbf{b}_{ks}^H \mathbf{b}_{ks} \le \bar{\check{p}}_{ks}, \forall k, s \quad (7)$$

where $P_{max}$, $\check{p}_n$, $\check{p}_k$ and $\bar{\check{p}}_{ks}$ are the maximum powers available at the BS, $n$th BS antenna, $k$th user and $k$th user $s$th symbol, respectively. Since, the problems $\mathcal{P}1$ - $\mathcal{P}4$ are not convex, convex optimization framework can not be applied to solve them. To the best of our knowledge, duality based solutions for the problems $\mathcal{P}1 - \mathcal{P}4$ are not known. In the following, we present duality based iterative algorithm to solve each of these problems. For better exposition of our algorithms, let us explain the general framework of downlink-uplink duality based approach for solving each of these problems as shown in **Algorithm I**.

**Algorithm I**

Initialization: For each problem, initialize $\{\mathbf{B}_k\}_{k=1}^K$ such that the power constraint functions are satisfied with equality. Then, update $\{\mathbf{W}_k\}_{k=1}^K$ by using minimum average-mean square-error (MAMSE) receiver approach, i.e.,

$$\mathbf{W}_k = (\boldsymbol{\Gamma}_k^{DL})^{-1} \widehat{\mathbf{H}}_k^H \mathbf{B}_k, \; \forall k. \quad (8)$$

**Repeat**
**Uplink channel**
1) Transfer the total sum AMSE from downlink to uplink channel.
2) Update the receivers of the uplink channel $\{\mathbf{T}_k\}_{k=1}^K$ using MAMSE receiver technique.

**Downlink channel**
3) Transfer the total sum AMSE from uplink to downlink channel.
4) Update the receivers of the downlink channel $\{\mathbf{W}_k\}_{k=1}^K$ by MAMSE receiver (8).

**Until** convergence.

One can notice that this iterative algorithm is already known in [3], [4] and [14]. However, the iterative approach of these papers are limited to solve MSE-based constrained with a total BS power problems with $\{\mathbf{R}_{nk} = \sigma^2 \mathbf{I}\}_{k=1}^K$. As will be clear later, the main challenge arises at step 3 of **Algorithm I** where at this step, the approaches of these papers are not able to ensure the power constraints of $\mathcal{P}1 - \mathcal{P}4$.

In the following sections, we establish our novel sum AMSE downlink-uplink and downlink-interference duality. The downlink-uplink duality are able to maintain the power constraints of $\mathcal{P}1$ and $\mathcal{P}2$ at step 3 of **Algorithm I**. Thus, the latter duality can solve these two problems using **Algorithm I**. The sum AMSE downlink-interference duality are able to maintain the power constraint of each user and symbol at step 3 of **Algorithm I**. Hence, these duality can solve the problems $\mathcal{P}3$ and $\mathcal{P}4$ with **Algorithm I**.

## V. Sum AMSE downlink-uplink duality to solve $\mathcal{P}1$

In this section, we establish the sum AMSE downlink-uplink duality to solve total BS power constrained robust sum MSE minimization problem ($\mathcal{P}1$). This duality can be established by assuming that $\mathbf{\Psi} = \sigma^2 \mathbf{I}$, where $\sigma^2 > 0$. Now, we transfer the sum AMSE from downlink to uplink channel and vice versa. To this end, we compute the total sum AMSE in the uplink channel as

$$\bar{\xi}^{UL} = S + \sum_{k=1}^{K} \text{tr}\{\mathbf{T}_k^H \mathbf{\Gamma}_c \mathbf{T}_k + \sigma^2 \mathbf{T}_k^H \mathbf{T}_k - 2\Re\{\mathbf{T}_k^H \hat{\mathbf{H}}_k \mathbf{V}_k\}\} \quad (9)$$

where $\mathbf{\Gamma}_c = \sum_{i=1}^{K} (\hat{\mathbf{H}}_i \mathbf{V}_i \mathbf{V}_i^H \hat{\mathbf{H}}_i^H + \sigma_{ei}^2 \text{tr}\{\mathbf{R}_{mi} \mathbf{V}_i \mathbf{V}_i^H\} \mathbf{R}_{bi})$.

### A. Sum AMSE transfer (From downlink to uplink channel)

The sum AMSE can be transferred from downlink to uplink channel by using a nonzero scaling factor $\tilde{\beta}$ which satisfies

$$\mathbf{V} = \tilde{\beta} \mathbf{W}, \quad \mathbf{T} = \mathbf{B}/\tilde{\beta} \quad (10)$$

where $\mathbf{V} = \text{blkdiag}(\mathbf{V}_1, \cdots, \mathbf{V}_K)$ and $\mathbf{T} = [\mathbf{T}_1, \cdots, \mathbf{T}_K]$. Substituting $\mathbf{V}$ and $\mathbf{T}$ of (9) by (10) and then equating $\bar{\xi}^{UL} = \bar{\xi}^{DL}$ yields

$$\text{tr}\{\mathbf{B}^H \hat{\mathbf{H}} \mathbf{W} \mathbf{W}^H \hat{\mathbf{H}}^H \mathbf{B}\} + \frac{1}{\tilde{\beta}^2} \sigma^2 \text{tr}\{\mathbf{B}^H \mathbf{B}\}$$
$$+ \sum_{k=1}^{K} \text{tr}\{\mathbf{B}_k^H (\sum_{i=1}^{K} \sigma_{ei}^2 \text{tr}\{\mathbf{R}_{mi} \mathbf{W}_i \mathbf{W}_i^H\} \mathbf{R}_{bi}) \mathbf{B}_k\} =$$
$$\text{tr}\{\mathbf{W}^H \hat{\mathbf{H}}^H \mathbf{B} \mathbf{B}^H \hat{\mathbf{H}} \mathbf{W}\} + \text{tr}\{\mathbf{W}^H \mathbf{R}_n \mathbf{W}\}$$
$$+ \sum_{k=1}^{K} \text{tr}\{\mathbf{W}_k^H (\sigma_{ek}^2 \text{tr}\{\mathbf{R}_{bk} \mathbf{B} \mathbf{B}^H\} \mathbf{R}_{mk}) \mathbf{W}_k\} \quad (11)$$

where $\mathbf{R}_n = \text{blkdiag}(\mathbf{R}_{n1}, \cdots, \mathbf{R}_{nK})$. It follows, $\tilde{\beta}$ can be determined as

$$\tilde{\beta}^2 = \sigma^2 \text{tr}\{\mathbf{B}^H \mathbf{B}\}/\text{tr}\{\mathbf{W}^H \mathbf{R}_n \mathbf{W}\}. \quad (12)$$

Thus, during the sum AMSE transfer from downlink to uplink channel, the following holds true

$$P_{sum}^{DL} = \text{tr}\{\mathbf{B}^H \mathbf{B}\} = \frac{1}{\sigma^2} \text{tr}\{\mathbf{V}^H \mathbf{R}_n \mathbf{V}\}. \quad (13)$$

### B. Sum AMSE transfer (From uplink to downlink channel)

For a given uplink channel precoder/decoder pairs ($\mathbf{V}/\mathbf{T}$), the sum AMSE can be transferred from uplink to downlink channel by using a positive $\beta$ which satisfies

$$\mathbf{B} = \beta \mathbf{T}, \quad \mathbf{W} = \mathbf{V}/\beta. \quad (14)$$

Substituting these $\mathbf{B}$ and $\mathbf{W}$ in (3) and equating $\xi^{UL} = \xi^{DL}$, $\beta$ can be determined as

$$\beta^2 = \text{tr}\{\mathbf{V}^H \mathbf{R}_n \mathbf{V}\}/\text{tr}\{\sigma^2 \mathbf{T}^H \mathbf{T}\}. \quad (15)$$

As can be seen from (12) and (15), the scaling factors $\tilde{\beta}$ and $\beta$ do not depend on $\{\sigma_{ek}^2\}_{k=1}^{K}$. This fact can be seen from $\bar{\xi}^{DL}$ and $\bar{\xi}^{UL}$, after substituting $\{\mathbf{\Gamma}_k^{DL}\}_{k=1}^{K}$ and $\mathbf{\Gamma}_c$, where $\{\sigma_{ek}^2\}_{k=1}^{K}$ are amplified by the same factor. The downlink power is given by

$$P_{sum}^{DL} = \text{tr}\{\mathbf{B} \mathbf{B}^H\}$$
$$= \text{tr}\{\beta^2 \mathbf{T} \mathbf{T}^H\} = \text{tr}\{\mathbf{V}^H \mathbf{R}_n \mathbf{V}\}/\sigma^2. \quad (16)$$

From (13) and (16), we can see that $\text{tr}\{\mathbf{B} \mathbf{B}^H\} = \text{tr}\{\mathbf{V}^H \mathbf{R}_n \mathbf{V}\}/\sigma^2 = P_{sum}^{DL}$ is always satisfied. Consequently, the duality of this section maintains the power constraint of $\mathcal{P}1$ during the sum AMSE transfer from uplink to downlink channel and vice versa. Therefore, with the duality of this section, the latter problem can be solved iteratively by **Algorithm I**. When $\{\sigma^2 = 1, \{\sigma_{ek}^2 = 0\}_{k=1}^{K}\}$, the duality of this section turns to that of in [15]. When $\mathbf{R}_n = \sigma^2 \mathbf{I}$, our duality fits that of in [13]. Hence, the sum AMSE duality of this section generalizes all existing sum AMSE downlink-uplink duality.

## VI. Sum AMSE downlink-uplink duality to solve $\mathcal{P}2$

This duality is established to solve per BS antenna power constrained robust sum MSE minimization problem ($\mathcal{P}2$). To establish this duality, we first compute the total sum AMSE in the uplink channel as

$$\bar{\xi}^{UL_2} = \text{tr}\{\mathbf{T}^H \mathbf{\Gamma}_c \mathbf{T} + \mathbf{T}^H \mathbf{\Psi} \mathbf{T} - \mathbf{T}^H \hat{\mathbf{H}} \mathbf{V} - \mathbf{V}^H \hat{\mathbf{H}}^H \mathbf{T}\} + S. \quad (17)$$

Now, we transfer the sum AMSE from downlink to uplink channel and vice versa.

### A. Sum AMSE transfer (From downlink to uplink channel)

In order to use this sum AMSE transfer for a per BS antenna power constrained robust sum MSE minimization problem, we set the uplink channel precoder and decoder pairs as

$$\mathbf{V} = \mathbf{W}, \quad \mathbf{T} = \mathbf{B}. \quad (18)$$

Substituting (18) into (17) and then equating $\bar{\xi}^{DL} = \bar{\xi}^{UL_2}$ yields

$$\sum_{k=1}^{K} \text{tr}\{\mathbf{W}_k^H \mathbf{R}_{nk} \mathbf{W}_k\} = \sum_{k=1}^{K} \text{tr}\{\mathbf{B}_k^H \mathbf{\Psi} \mathbf{B}_k\},$$
$$\Rightarrow \tau = \sum_{n=1}^{N} \psi_n \tilde{p}_n = \tilde{\mathbf{p}}^T \boldsymbol{\psi} \quad (19)$$

where $\tau = \sum_{k=1}^{K} \text{tr}\{\mathbf{W}_k^H \mathbf{R}_{nk} \mathbf{W}_k\}$, $\boldsymbol{\psi} = [\psi_1, \cdots, \psi_N]^T$, $\tilde{\mathbf{p}} = [\tilde{p}_1, \cdots, \tilde{p}_N]^T$, $\tilde{p}_n = \text{tr}\{\tilde{\mathbf{b}}_n^H \tilde{\mathbf{b}}_n\}$ and $\tilde{\mathbf{b}}_n^H$ is the $n$th row of $\mathbf{B}$. The above equation shows that by choosing $\{\psi_n\}_{n=1}^N$ appropriately, we can transfer any precoder/decoder pairs of the downlink channel to the corresponding decoder/precoder pairs of the uplink channel. However, here $\{\psi_n\}_{n=1}^N$ should be selected in a way that $\mathcal{P}2$ can be solved using **Algorithm I**. To this end, we choose $\boldsymbol{\psi}$ as

$$\tau \geq \tilde{\mathbf{p}}^T \boldsymbol{\psi}. \quad (20)$$

By doing so, the uplink channel will achieve lower sum AMSE compared to that of the downlink channel ($\bar{\xi}^{DL} \geq \bar{\xi}^{UL_2}$). As will be clear later, to solve (5) with **Algorithm I**, $\boldsymbol{\psi}$ should be selected ensuring (20). This shows that for $\mathcal{P}2$, step 1 of **Algorithm I** can be carried out with (18).

To perform step 2 of **Algorithm I**, we update $\mathbf{T}$ of (18) by using the uplink MAMSE receiver which can be expressed as

$$\mathbf{T} = (\boldsymbol{\Gamma}_c + \boldsymbol{\Psi})^{-1} \widehat{\mathbf{H}} \mathbf{V} = (\mathbf{A} + \boldsymbol{\Upsilon} + \boldsymbol{\Psi})^{-1} \widehat{\mathbf{H}} \mathbf{W} \quad (21)$$

where $\mathbf{A} = \widehat{\mathbf{H}} \mathbf{W} \mathbf{W}^H \widehat{\mathbf{H}}^H$, $\boldsymbol{\Upsilon} = \sum_{i=1}^{K} \sigma_{ei}^2 \text{tr}\{\mathbf{R}_{mi} \mathbf{W}_i \mathbf{W}_i^H\} \mathbf{R}_{bi}$ and the second equality is obtained by applying (18) (i.e., $\mathbf{V}=\mathbf{W}$). The above expression shows that by choosing $\{\psi_n > 0\}_{n=1}^N$, $(\mathbf{A} + \boldsymbol{\Upsilon} + \boldsymbol{\Psi})$ is always invertible. Next, we transfer the total sum AMSE from uplink to downlink channel (i.e., we perform step 3 of **Algorithm I**) and show that the latter sum AMSE transfer ensures the power constraint of each BS antenna.

### B. Sum AMSE transfer (From uplink to downlink channel)

For a given total sum AMSE in the uplink channel, we can achieve the same sum AMSE in the downlink channel by using a nonzero scaling factor ($\breve{\beta}$) which satisfies

$$\widetilde{\mathbf{B}} = \breve{\beta} \mathbf{T}, \quad \widetilde{\mathbf{W}} = \mathbf{V}/\breve{\beta}. \quad (22)$$

In this precoder/decoder transformation, we use the notations $\widetilde{\mathbf{B}}$ and $\widetilde{\mathbf{W}}$ to differentiate with the precoder and decoder matrices used in Section VI-A. By substituting (22) into (3) (with $\widetilde{\mathbf{B}}=\mathbf{B}$, $\widetilde{\mathbf{W}}=\mathbf{W}$) and then equating the resulting sum AMSE with that of the uplink channel (17), $\breve{\beta}$ can be determined as

$$\frac{1}{\breve{\beta}^2} \sum_{k=1}^{K} \text{tr}\{\mathbf{V}_k^H \mathbf{R}_{nk} \mathbf{V}_k\} = \sum_{n=1}^{N} \psi_n \tilde{\mathbf{t}}_n^H \tilde{\mathbf{t}}_n,$$

$$\Rightarrow \breve{\beta}^2 = \sum_{k=1}^{K} \text{tr}\{\mathbf{W}_k^H \mathbf{R}_{nk} \mathbf{W}_k\} / \sum_{n=1}^{N} \psi_n \tilde{\mathbf{t}}_n^H \tilde{\mathbf{t}}_n$$

$$= \frac{\tau}{\sum_{n=1}^{N} \psi_n \tilde{\mathbf{t}}_n^H \tilde{\mathbf{t}}_n} \quad (23)$$

where $\tilde{\mathbf{t}}_n^H$ is the $n$th row of the MAMSE matrix $\mathbf{T}$ (21) which is given by $[(\mathbf{A}+\boldsymbol{\Upsilon}+\boldsymbol{\Psi})^{-1}]_{(n,:)} \widehat{\mathbf{H}} \mathbf{W}$ and the second equality follows from (18).

The power of each BS antenna in the downlink channel is thus given by

$$\tilde{\mathbf{b}}_n^H \tilde{\mathbf{b}}_n = \text{tr}\{\breve{\beta}^2 \tilde{\mathbf{t}}_n^H \tilde{\mathbf{t}}_n\} = \frac{\tau \tilde{\mathbf{t}}_n^H \tilde{\mathbf{t}}_n}{\sum_{i=1}^{N} \psi_i \tilde{\mathbf{t}}_i^H \tilde{\mathbf{t}}_i} \leq \breve{p}_n, \forall n \quad (24)$$

where $\tilde{\mathbf{b}}_n^H$ is the $n$th row of $\widetilde{\mathbf{B}}$. We can rewrite the above expression as

$$\psi_n \geq f_n, \forall n \quad (25)$$

where

$$f_n = \frac{\tau}{\breve{p}_n} \frac{\psi_n \tilde{\mathbf{t}}_n^H \tilde{\mathbf{t}}_n}{\sum_{i=1}^{N} \psi_i \tilde{\mathbf{t}}_i^H \tilde{\mathbf{t}}_i} = \frac{\tau}{\breve{p}_n} \times$$

$$\frac{\psi_n [(\mathbf{A}+\boldsymbol{\Upsilon}+\boldsymbol{\Psi})^{-1}]_{(n,:)} \mathbf{A} ([(\mathbf{A}+\boldsymbol{\Upsilon}+\boldsymbol{\Psi})^{-1}]_{(n,:)})^H}{\sum_{i=1}^{N} \psi_i [(\mathbf{A}+\boldsymbol{\Upsilon}+\boldsymbol{\Psi})^{-1}]_{(i,:)} \mathbf{A} ([(\mathbf{A}+\boldsymbol{\Upsilon}+\boldsymbol{\Psi})^{-1}]_{(i,:)})^H}.$$

Let us assume that there exist $\{\psi_n > 0\}_{n=1}^N$ that satisfy

$$\psi_n = f_n, \quad \forall n. \quad (26)$$

From (26), we get

$$\psi_n = \frac{\tau}{\breve{p}_n} \frac{\psi_n \tilde{\mathbf{t}}_n^H \tilde{\mathbf{t}}_n}{\sum_{i=1}^{N} \psi_i \tilde{\mathbf{t}}_i^H \tilde{\mathbf{t}}_i} \Rightarrow \psi_n \breve{p}_n = \tau \frac{\psi_n \tilde{\mathbf{t}}_n^H \tilde{\mathbf{t}}_n}{\sum_{i=1}^{N} \psi_i \tilde{\mathbf{t}}_i^H \tilde{\mathbf{t}}_i}$$

$$\Rightarrow \sum_{n=1}^{N} \psi_n \breve{p}_n = \sum_{n=1}^{N} \tau \frac{\psi_n \tilde{\mathbf{t}}_n^H \tilde{\mathbf{t}}_n}{\sum_{i=1}^{N} \psi_i \tilde{\mathbf{t}}_i^H \tilde{\mathbf{t}}_i} = \tau. \quad (27)$$

The above expression shows that the solution of (26) satisfies (27). Moreover, as $\{\breve{p}_n \geq \tilde{p}_n\}_{n=1}^N$, the latter solution also satisfies (20). Therefore, for $\mathcal{P}2$, by choosing $\{\psi_n\}_{n=1}^N$ such that (26) is satisfied, step 3 of **Algorithm I** can be performed. Next, we show that there exists at least a set of feasible $\{\psi_n > 0\}_{n=1}^N$ that satisfy (26). In this regard, we consider the following Theorem [23], [24].

*Theorem 1:* Let $(\mathbf{X}, \|.\|_2)$ be a complete metric space. We say that $F : \mathbf{X} \to \mathbf{X}$ is an almost contraction, if there exist $\kappa \in [0, 1)$ and $\chi \geq 0$ such that

$$\|F(\mathbf{x}) - F(\mathbf{y})\|_2 \leq \kappa \|\mathbf{x} - \mathbf{y}\|_2 + \chi \|\mathbf{y} - F(\mathbf{x})\|_2,$$
$$\forall \mathbf{x}, \mathbf{y} \in \mathbf{X}. \quad (28)$$

If $F$ satisfies (28), then the following holds true
1) $\exists \mathbf{x} \in \mathbf{X} : \mathbf{x} = F(\mathbf{x})$.
2) For any initial $\mathbf{x}_0 \in \mathbf{X}$, the iteration $\mathbf{x}_{n+1} = F(\mathbf{x}_n)$ for n = 0, 1, 2, $\cdots$ converges to some $\mathbf{x}^\star \in \mathbf{X}$.
3) The solution $\mathbf{x}^\star$ is not necessarily unique.

*Proof:* See [23], [24] (Theorem 1.1).

Note that in [23] and [24], *Theorem 1* has been proven for a generalized complete metric space $(\mathbf{X}, d)$ instead of $(\mathbf{X}, \|.\|_2)$. By defining $F$ as $F(\boldsymbol{\psi}) \triangleq [f_1, f_2, \cdots, f_N]$ with $\{\psi_n \in [\epsilon, (\tau - \epsilon \sum_{i=1, i \neq n}^N \breve{p}_i)/\breve{p}_n]\}_{n=1}^N$[3], it can be easily seen that $\|F(\boldsymbol{\psi}_1) - F(\boldsymbol{\psi}_2)\|_2$, $\|\boldsymbol{\psi}_1 - \boldsymbol{\psi}_2\|_2$ and $\|\boldsymbol{\psi}_2 - F(\boldsymbol{\psi}_1)\|_2$ are bounded for any $\boldsymbol{\psi}_1, \boldsymbol{\psi}_2 \in \boldsymbol{\psi}$. This shows the existence of $\kappa$ and $\chi$ ensuring (28). Consequently, $F(\boldsymbol{\psi})$ is an almost contraction which implies

$$\boldsymbol{\psi}_{n+1} = F(\boldsymbol{\psi}_n), \text{ with } \boldsymbol{\psi}_0 = [\psi_{01}, \psi_{02}, \cdots, \psi_{0N}] \geq \epsilon \mathbf{1}_N,$$
$$\text{for } n = 0, 1, 2, \cdots \text{ converges} \quad (29)$$

where $\mathbf{1}_N$ is an $N$ length vector with each element equal to unity. Thus, there exists $\{\psi_n \geq \epsilon\}_{n=1}^N$ that satisfy (26) and its solution can be obtained using the above fixed point iterations. Once the appropriate $\{\psi_n\}_{n=1}^N$ is obtained, step 4 of **Algorithm I** is immediate and hence $\mathcal{P}2$ can be solved iteratively using this algorithm. Note that when $\{\sigma_{ek}^2 = 0\}_{k=1}^K$, the duality of this section turns to that of [20].

---
[3]Here $\epsilon$ should be very close to zero (for our simulation, we use $\epsilon = \min(10^{-6}, \{\tau/\breve{p}_n\}_{n=1}^N)$).

## VII. SUM AMSE DOWNLINK-INTERFERENCE DUALITY TO SOLVE $\mathcal{P}3$

This kind of duality is established to solve per user power constrained robust sum MSE minimization problem ($\mathcal{P}3$). To establish this duality, it is sufficient to assume $\{\mathbf{\Delta}_{ks} = \mathbf{\Delta}_k = \mu_k \mathbf{I}_N, \forall s\}_{k=1}^K$. With this assumption, the $k$th user AMSE and total sum AMSE in the interference channel can be expressed as

$$\bar{\xi}_k^{I_1} = \mathrm{tr}\{\mathbf{T}_k^H \mathbf{\Gamma}_c \mathbf{T}_k + \mu_k \mathbf{T}_k^H \mathbf{T}_k + \mathbf{I}_{S_k} - \mathbf{T}_k^H \mathbf{H}_k \mathbf{V}_k - \mathbf{V}_k^H \mathbf{H}_k^H \mathbf{T}_k\} \tag{30}$$

$$\bar{\xi}^{I_1} = \sum_{k=1}^K \bar{\xi}_k^{I_1} = \mathrm{tr}\{\mathbf{T}^H \mathbf{\Gamma}_c \mathbf{T} - \mathbf{T}^H \mathbf{H} \mathbf{V} - \mathbf{V}^H \mathbf{H}^H \mathbf{T}\} + S + \sum_{k=1}^K \mathrm{tr}\{\mu_k \mathbf{T}_k^H \mathbf{T}_k\}. \tag{31}$$

Like in Section VI, here we transfer the sum AMSE from downlink to interference channel and vice versa.

### A. Sum AMSE transfer (From downlink to interference channel)

Like in Section VI, it can be shown that the sum AMSE can be transferred from downlink to interference channel by setting $\mathbf{V} = \mathbf{W}$ and $\mathbf{T} = \mathbf{B}$. Substituting these $\mathbf{V}$ and $\mathbf{T}$ in (31), then equating $\bar{\xi}^{DL} = \bar{\xi}^{I_1}$ yields

$$\sum_{k=1}^K \mathrm{tr}\{\mathbf{W}_k^H \mathbf{R}_{nk} \mathbf{W}_k\} = \sum_{k=1}^K \mathrm{tr}\{\mu_k \mathbf{B}_k^H \mathbf{B}_k\} \Rightarrow \tau = \tilde{\mathbf{p}}^T \boldsymbol{\mu}$$

where $\tilde{p}_k = \mathrm{tr}\{\mathbf{B}_k^H \mathbf{B}_k\}$, $\tilde{\mathbf{p}} = [\tilde{p}_1, \cdots, \tilde{p}_K]^T$ and $\boldsymbol{\mu} = [\mu_1, \cdots, \mu_K]^T$. For this sum AMSE transfer, we also choose $\boldsymbol{\mu}$ such that

$$\tau \geq \tilde{\mathbf{p}}^T \boldsymbol{\mu} \tag{32}$$

is satisfied. This shows that step 1 of **Algorithm I** can be carried out by choosing $\mathbf{V} = \mathbf{W}$ and $\mathbf{T} = \mathbf{B}$.

To perform step 2 of **Algorithm I**, we update each receiver of the interference channel by applying MAMSE receiver method which is given by

$$\mathbf{t}_{ks} = (\mathbf{\Gamma}_c + \mu_k \mathbf{I}_N)^{-1} \widehat{\mathbf{H}}_k \mathbf{v}_{ks}$$
$$= (\mathbf{A} + \mathbf{\Upsilon} + \mu_k \mathbf{I}_N)^{-1} \widehat{\mathbf{H}}_k \mathbf{w}_{ks}, \ \forall k, s. \tag{33}$$

The above expression shows that by choosing $\{\mu_k > 0\}_{k=1}^K$, invertibility of $(\mathbf{A} + \mathbf{\Upsilon} + \mu_k \mathbf{I})$ can be ensured.

### B. Sum AMSE transfer (From interference to downlink channel)

For a given total sum AMSE in the interference channel, we can achieve the same sum AMSE in the downlink channel by using a nonzero scaling factor ($\check{\beta}$) which satisfies $\widetilde{\mathbf{B}} = \check{\beta}\mathbf{T}, \widetilde{\mathbf{W}} = \mathbf{V}/\check{\beta}$. By substituting $\widetilde{\mathbf{B}}$ and $\widetilde{\mathbf{W}}$ in (3) (with $\widetilde{\mathbf{B}}=\mathbf{B}$, $\widetilde{\mathbf{W}}=\mathbf{W}$) and then equating the resulting total sum AMSE with that of the interference channel (31), $\check{\beta}$ can be determined as

$$\check{\beta}^2 = \sum_{k=1}^K \mathrm{tr}\{\mathbf{W}_k^H \mathbf{R}_{nk} \mathbf{W}_k\} / \sum_{k=1}^K \mathrm{tr}\{\mu_k \mathbf{T}_k^H \mathbf{T}_k\}$$
$$= \frac{\tau}{\sum_{k=1}^K \mu_k \mathrm{tr}\{(\mathbf{A} + \mathbf{\Upsilon} + \mu_k \mathbf{I}_N)^{-2} \mathbf{A}_k\}} \tag{34}$$

where $\mathbf{A}_k = \widehat{\mathbf{H}}_k \mathbf{W}_k \mathbf{W}_k^H \widehat{\mathbf{H}}_k^H$ and the second equality is derived by substituting $\{\mathbf{t}_{ks}, \forall s\}_{k=1}^K$ with the MAMSE receiver (33). The power of each user in the downlink channel is thus given by

$$\mathrm{tr}\{\widetilde{\mathbf{B}}_k \widetilde{\mathbf{B}}_k^H\} = \mathrm{tr}\{\check{\beta}^2 \mathbf{T}_k \mathbf{T}_k^H\} \tag{35}$$
$$= \frac{\tau \mathrm{tr}\{(\mathbf{A} + \mathbf{\Upsilon} + \mu_k \mathbf{I}_N)^{-2} \mathbf{A}_k\}}{\sum_{i=1}^K \mu_i \mathrm{tr}\{(\mathbf{A} + \mathbf{\Upsilon} + \mu_i \mathbf{I}_N)^{-2} \mathbf{A}_i\}} \leq \check{p}_k.$$

This inequality constraint can be expressed as

$$\mu_k \geq \check{f}_k, \ \forall k \tag{36}$$

where $\check{f}_k = \frac{\tau}{\check{p}_k} \frac{\mathrm{tr}\{\mu_k(\mathbf{A}+\mathbf{\Upsilon}+\mu_k\mathbf{I}_N)^{-2}\mathbf{A}_k\}}{\sum_{i=1}^K \mu_i \mathrm{tr}\{(\mathbf{A}+\mathbf{\Upsilon}+\mu_i\mathbf{I}_N)^{-2}\mathbf{A}_i\}}$. Like in Section VI, it can be shown that there exist $\{\mu_k > 0\}_{k=1}^K$ that satisfy

$$\mu_k = \check{f}_k, \ \forall k. \tag{37}$$

Moreover, (37) can be solved exactly like that of (26) and the solution of (37) satisfies (32). As a result, one can apply the duality of this section to solve $\mathcal{P}3$ by **Algorithm I**.

## VIII. SUM AMSE DOWNLINK-INTERFERENCE DUALITY TO SOLVE $\mathcal{P}4$

This duality is established to solve per symbol power constrained robust sum MSE minimization problem ($\mathcal{P}4$). To establish this duality, it is sufficient to assume $\{\mathbf{\Delta}_{ks} = \tilde{\mu}_{ks}\mathbf{I}_N, \forall s\}_{k=1}^K$. With the latter assumption, the $k$th user $s$th symbol AMSE and total sum AMSE in the interference channel can be expressed as

$$\bar{\xi}_{ks}^{I_2} = 1 + \mathbf{t}_{ks}^H \mathbf{\Gamma}_c \mathbf{t}_{ks} + \tilde{\mu}_{ks} \mathbf{t}_{ks}^H \mathbf{t}_{ks} - 2\Re\{\mathbf{t}_{ks}^H \mathbf{H}_k \mathbf{v}_{ks}\} \tag{38}$$

$$\bar{\xi}^{I_2} = \sum_{k=1}^K \sum_{s=1}^{S_k} \bar{\xi}_{ks}^{I_2} = \mathrm{tr}\{\mathbf{T}^H \mathbf{\Gamma}_c \mathbf{T} - \mathbf{T}^H \mathbf{H} \mathbf{V} - \mathbf{V}^H \mathbf{H}^H \mathbf{T}\} + S + \sum_{k=1}^K \sum_{s=1}^{S_k} \tilde{\mu}_{ks} \mathbf{t}_{ks}^H \mathbf{t}_{ks}. \tag{39}$$

Like in Section VII, we now transfer the sum AMSE from downlink to interference channel and vice versa.

## A. Sum AMSE transfer (From downlink to interference channel)

The sum AMSE can be transferred from downlink to interference channel by setting $\mathbf{V} = \mathbf{W}$ and $\mathbf{T} = \mathbf{B}$. Substituting these $\mathbf{V}$ and $\mathbf{T}$ in (39), equating $\bar{\xi}^{DL} = \bar{\xi}^{I_2}$ and after some steps we get $\tau = \check{\mathbf{p}}^T \tilde{\boldsymbol{\mu}}$, where $\check{p}_{ks} = \mathbf{b}_{ks}^H \mathbf{b}_{ks}$, $\check{\mathbf{p}} = [\check{p}_{11}, \cdots, \check{p}_{1S_1}, \cdots \check{p}_{K1}, \cdots, \check{p}_{KS_K}]^T$ and $\tilde{\boldsymbol{\mu}} = [\tilde{\mu}_{11}, \cdots, \tilde{\mu}_{1S_1}, \cdots, \tilde{\mu}_{K1}, \cdots, \tilde{\mu}_{KS_K}]^T$. For this sum AMSE transfer, we also choose $\tilde{\boldsymbol{\mu}}$ such that

$$\tau \geq \check{\mathbf{p}}^T \tilde{\boldsymbol{\mu}} \qquad (40)$$

is satisfied. At this stage step 1 of **Algorithm I** is performed.

To reduce the sum AMSE in the interference channel (i.e., to perform step 2 of **Algorithm I**), we apply the MAMSE receiver for the $k$th user $s$th symbol as

$$\begin{aligned}\mathbf{t}_{ks} &= (\boldsymbol{\Gamma}_c + \tilde{\mu}_{ks}\mathbf{I}_N)^{-1}\widehat{\mathbf{H}}_k \mathbf{v}_{ks} \\ &= (\mathbf{A} + \boldsymbol{\Upsilon} + \tilde{\mu}_{ks}\mathbf{I}_N)^{-1}\widehat{\mathbf{H}}_k \mathbf{w}_{ks}, \ \forall k, s.\end{aligned} \qquad (41)$$

As we can see, when $\{\tilde{\mu}_{ks} > 0, \forall s\}_{k=1}^K$, $(\mathbf{A} + \boldsymbol{\Upsilon} + \tilde{\mu}_{ks}\mathbf{I}_N)$ is always invertible.

## B. Sum AMSE transfer (From interference to downlink channel)

For a given total sum AMSE in the interference channel, we can achieve the same sum AMSE in the downlink channel by using a nonzero scaling factor $\bar{\bar{\beta}}$ which satisfies $\widetilde{\mathbf{B}} = \bar{\bar{\beta}}\mathbf{T}$ and $\widetilde{\mathbf{W}} = \mathbf{V}/\bar{\bar{\beta}}$. By substituting these $\widetilde{\mathbf{B}}$ and $\widetilde{\mathbf{W}}$ in (3) (with $\widetilde{\mathbf{B}}=\mathbf{B}$, $\widetilde{\mathbf{W}}=\mathbf{W}$) and then equating the resulting total sum AMSE with that of the interference channel (39), $\bar{\bar{\beta}}$ can be determined as

$$\begin{aligned}\bar{\bar{\beta}}^2 &= \sum_{k=1}^K \mathrm{tr}\{\mathbf{W}_k^H \mathbf{R}_{nk} \mathbf{W}_k\} / \underbrace{\sum_{k=1}^K \sum_{s=1}^{S_k} \mathrm{tr}\{\tilde{\mu}_{ks} \mathbf{t}_{ks}^H \mathbf{t}_{ks}\}}_{\tau} \\ &= \frac{\tau}{\sum_{k=1}^K \sum_{s=1}^{S_k} \tilde{\mu}_{ks}\mathrm{tr}\{(\mathbf{A} + \boldsymbol{\Upsilon} + \tilde{\mu}_{ks}\mathbf{I}_N)^{-2} \mathbf{A}_{ks}\}}\end{aligned} \qquad (42)$$

where the second equality is derived by substituting $\{\mathbf{t}_{ks}, \forall s\}_{k=1}^K$ with the MAMSE receiver (41) and $\mathbf{A}_{ks} = \widehat{\mathbf{H}}_k \mathbf{w}_{ks} \mathbf{w}_{ks}^H \widehat{\mathbf{H}}_k^H$. The power of the $k$th user $s$th symbol in the downlink channel is thus given by

$$\mathrm{tr}\{\widetilde{\mathbf{b}}_{ks} \widetilde{\mathbf{b}}_{ks}^H\} = \mathrm{tr}\{\bar{\bar{\beta}}^2 \mathbf{t}_{ks} \mathbf{t}_{ks}^H\} \qquad (43)$$
$$= \frac{\tau \mathrm{tr}\{(\mathbf{A} + \boldsymbol{\Upsilon} + \tilde{\mu}_{ks}\mathbf{I}_N)^{-2} \mathbf{A}_{ks}\}}{\sum_{i=1}^K \sum_{j=1}^{S_i} \tilde{\mu}_{ij} \mathrm{tr}\{(\mathbf{A} + \boldsymbol{\Upsilon} + \tilde{\mu}_{ij}\mathbf{I}_N)^{-2} \mathbf{A}_{ij}\}} \leq \check{\bar{p}}_{ks}.$$

This inequality constraint can also be expressed as

$$\tilde{\mu}_{ks} \geq \check{\tilde{f}}_{ks}, \ \forall k, s \qquad (44)$$

where $\check{\tilde{f}}_{ks} = \frac{\tau \mathrm{tr}\{\tilde{\mu}_{ks}(\mathbf{A}+\boldsymbol{\Upsilon}+\tilde{\mu}_{ks}\mathbf{I}_N)^{-2}\mathbf{A}_{ks}\}}{\check{\bar{p}}_{ks}\sum_{i=1}^K\sum_{j=1}^{S_i}\tilde{\mu}_{ij}\mathrm{tr}\{(\mathbf{A}+\boldsymbol{\Upsilon}+\tilde{\mu}_{ij}\mathbf{I}_N)^{-2}\mathbf{A}_{ij}\}}$. Like in Section VI, it can be shown that there exist $\{\tilde{\mu}_{ks} > 0, \forall s\}_{k=1}^K$ that satisfy $\tilde{\mu}_{ks} = \check{\tilde{f}}_{ks}, \forall k, s$. Thus, the feasible $\{\tilde{\mu}_{ks}, \forall s\}_{k=1}^K$ of (44) can be found exactly like that of (26). As a result, with the duality of this section, we can solve $\mathcal{P}4$ iteratively with **Algorithm I**.

## IX. A METHOD TO IMPROVE THE CONVERGENCE SPEED OF **Algorithm I**

To increase the convergence speed of **Algorithm I**, a downlink power allocation step can be included inside **Algorithm I**. In [13], for fixed transmit and receive filters, the power allocation part of $\mathcal{P}1$ with $\{\mathbf{R}_{nk} = \sigma^2 \mathbf{I}, \widetilde{\mathbf{R}}_{mk} = \mathbf{I}\}_{k=1}^K$ has been formulated as a GP by employing the system model of [4]. This GP formulation is derived by extending the approach of [4] to imperfect CSI scenario. Moreover, in [14], we show that the system model of [4] is appropriate to solve any kind of total BS power constrained robust MSE-based problems using duality approach (alternating optimization). This motivates us to utilize the system model of [4] in the downlink channel only and then include the power allocation step (i.e., GP) into **Algorithm I** for each of the problems $\mathcal{P}1 - \mathcal{P}4$. To this end, we decompose the precoders and decoders of the downlink channel as

$$\mathbf{B}_k = \mathbf{G}_k \mathbf{P}_k^{1/2}, \quad \mathbf{W}_k = \mathbf{U}_k \boldsymbol{\alpha}_k \mathbf{P}_k^{-1/2}, \ \forall k \qquad (45)$$

where $\mathbf{P}_k = \mathrm{diag}(p_{k1}, \cdots, p_{kS_k}) \in \Re^{S_k \times S_k}$, $\mathbf{G}_k = [\mathbf{g}_{k1} \cdots \mathbf{g}_{kS_k}] \in \mathcal{C}^{N \times S_k}$, $\mathbf{U}_k = [\mathbf{u}_{k1} \cdots \mathbf{u}_{kS_k}] \in \mathcal{C}^{M_k \times S_k}$ and $\boldsymbol{\alpha}_k = \mathrm{diag}(\alpha_{k1}, \cdots, \alpha_{kS_k}) \in \Re^{S_k \times S_k}$ are the transmit power, unity norm transmit filter, unity norm receive filter and receiver scaling factor matrices of the $k$th user, respectively, i.e., $\{\mathbf{g}_{ks}^H \mathbf{g}_{ks} = \mathbf{u}_{ks}^H \mathbf{u}_{ks} = 1, \forall s\}_{k=1}^K$. With this decomposition, our downlink system model (Fig. 1.(a)) can be plotted as shown in Fig. 2.

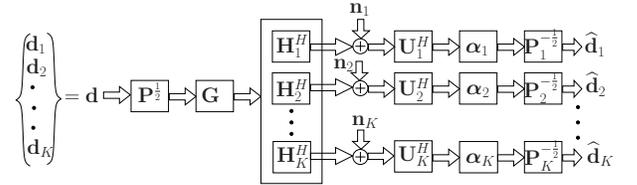

Fig. 2. Downlink multiuser MIMO system model 2.

For better exposition of the GP, we collect the powers of the entire symbol as $\mathbf{P} = \mathrm{blkdiag}(\mathbf{P}_1, \cdots, \mathbf{P}_K) = \mathrm{diag}(p_1, \cdots, p_S)$, the overall filter matrix at the BS and users as $\mathbf{G} = [\mathbf{G}_1, \cdots, \mathbf{G}_K] = [\mathbf{g}_1, \cdots, \mathbf{g}_S]$ and $\mathbf{U} = \mathrm{blkdiag}(\mathbf{U}_1, \cdots, \mathbf{U}_K) = [\mathbf{u}_1, \cdots, \mathbf{u}_S]$, respectively, where $\mathbf{g}_l \in C^{N \times 1} (\mathbf{u}_l \in \mathcal{C}^{M \times 1})$ is the transmit (receive) filter of the $l$th symbol with $\|\mathbf{g}_l\|_2 = \|\mathbf{u}_l\|_2 = 1$. The scaling factors are stacked as $\boldsymbol{\alpha} = \mathrm{blkdiag}(\boldsymbol{\alpha}_1, \cdots, \boldsymbol{\alpha}_K) = \mathrm{diag}(\alpha_1, \cdots, \alpha_S)$ and the AMSE in the downlink channel can be collected as $\bar{\boldsymbol{\xi}} = [\bar{\xi}_{1,1}^{DL}, \cdots, \bar{\xi}_{K,S_K}^{DL}]^T = [\bar{\xi}_1^{DL}, \cdots, \bar{\xi}_S^{DL}]^T = [\{\bar{\xi}_l^{DL}\}_{l=1}^S]^T$ (refer [4] and [13] for the details about (45) and the above descriptions). By applying (45) and the technique of [4] (see (23) of [4]), the AMSE of the $l$th symbol in the downlink channel can be expressed as

$$\bar{\xi}_l^{DL} = p_l^{-1}[(\mathbf{D} + \boldsymbol{\alpha}^2 \boldsymbol{\Phi}^T)\mathbf{p}]_l + p_l^{-1}\alpha_l^2 \mathbf{u}_l^H \mathbf{R}_n \mathbf{u}_l \qquad (46)$$

where

$$[\mathbf{\Phi}]_{l,j} = \begin{cases} \sigma_{ef(j)}^2 \|\mathbf{R}_{mf(j)}^{1/2} \mathbf{u}_j\|_2^2 \|\mathbf{R}_{bf(j)}^{1/2} \mathbf{g}_l\|_2^2 \\ + |\mathbf{g}_l^H \widehat{\mathbf{H}} \mathbf{u}_j|^2, & \text{for } l \neq j \\ 0, & \text{for } l = j \end{cases} \quad (47)$$

$$[\mathbf{D}]_{l,l} = \alpha_l^2 (|\mathbf{g}_l^H \widehat{\mathbf{H}} \mathbf{u}_l|^2 + \sigma_{ef(l)}^2 \|\mathbf{R}_{mf(l)}^{1/2} \mathbf{u}_l\|_2^2 \|\mathbf{R}_{bf(l)}^{1/2} \mathbf{g}_l\|_2^2) - 2\alpha_l \Re(\mathbf{u}_l^H \widehat{\mathbf{H}}^H \mathbf{g}_l) + 1 \quad (48)$$

and $\mathbf{p} = [p_1, \cdots, p_S]^T$. In (47) and (48), $f(i)$ is the smallest $k$ such that $\sum_{m=1}^k S_m - i \geq 0$, where $S_m$ is the $m$th user total number of symbols. For example, if $K = 2$, $M_k = 2$ and $S_k = 2$ then $\sigma_{ef(1)}^2 = \sigma_{ef(2)}^2 = \sigma_{e1}^2$, $\sigma_{ef(3)}^2 = \sigma_{ef(4)}^2 = \sigma_{e2}^2$, $\mathbf{R}_{bf(1)} = \mathbf{R}_{bf(2)} = \mathbf{R}_{b1}$, $\mathbf{R}_{mf(1)} = \mathbf{R}_{mf(2)} = \mathbf{R}_{m1}$, $\mathbf{R}_{bf(3)} = \mathbf{R}_{bf(4)} = \mathbf{R}_{b2}$ and $\mathbf{R}_{mf(3)} = \mathbf{R}_{mf(4)} = \mathbf{R}_{m2}$. Using the above $\bar{\xi}_l^{DL}$, for fixed $\mathbf{G}, \mathbf{U}$ and $\boldsymbol{\alpha}$, the power allocation part of problem $\mathcal{P}1$ can be formulated as

$$\min_{\{p_l\}_{l=1}^S} \sum_{l=1}^S \bar{\xi}_l^{DL}, \quad \text{s.t} \sum_{l=1}^S p_l \leq P_{max}. \quad (49)$$

As $\bar{\xi}_l^{DL}$ is a posynomial (where $\{p_l\}_{l=1}^S$ are the variables), (49) is a GP for which global optimality is guaranteed. Thus, it can be efficiently solved by using standard interior point methods with a worst-case polynomial-time complexity [25]. Like in $\mathcal{P}1$, it can be shown that for fixed $\mathbf{G}, \mathbf{U}$ and $\boldsymbol{\alpha}$, the power allocation parts of $\mathcal{P}2 - \mathcal{P}4$ can be formulated as GPs. Our duality based iterative algorithm for each of these problems including the power allocation step is summarized in **Algorithm II**.

> **Algorithm II**
> Initialization: For $\mathcal{P}1$, set any $\sigma^2 > 0$ (we use $\sigma^2 = 1$ for our simulation). Then, initialize all the other parameters like in **Algorithm I**.
> **Repeat**
> **Uplink (Interference) channel**
> 1) For $\mathcal{P}1$ transfer the total sum AMSE from downlink to uplink channel by (12) and for $\mathcal{P}2 - \mathcal{P}4$ set $\mathbf{V} = \mathbf{W}, \mathbf{T} = \mathbf{B}$. Then, for $\mathcal{P}2$ compute $\{\psi_n\}_{n=1}^N$ using (25); for $\mathcal{P}3$ compute $\{\mu_k\}_{k=1}^K$ using (36) and for $\mathcal{P}4$ compute $\{\tilde{\mu}_{ks}, \forall s\}_{k=1}^K$ using (44).
> 2) Update the MAMSE receivers of the uplink (interference) channel of $\mathcal{P}1$, $\mathcal{P}2$, $\mathcal{P}3$ and $\mathcal{P}4$ using $\mathbf{T} = \boldsymbol{\Gamma}_c^{-1} \widehat{\mathbf{H}} \mathbf{V}$, (21), (33) and (41), respectively.
> **Downlink channel**
> 3) Transfer the total sum AMSE from uplink (interference) to downlink channel using (15), (23), (34) and (42) for $\mathcal{P}1$, $\mathcal{P}2$, $\mathcal{P}3$ and $\mathcal{P}4$ respectively.
> 4) For each of the problems $\mathcal{P}1 - \mathcal{P}4$, decompose the precoder and decoder matrices of each user as in (45). Then, formulate and solve the GP power allocation part. For example, the power allocation part of $\mathcal{P}1$ can be expressed in GP form as (49).
> 5) For each of the problems $\mathcal{P}1 - \mathcal{P}4$, by keeping $\{\mathbf{P}_k\}_{k=1}^K$ constant, update the receive filters $\{\mathbf{U}_k\}_{k=1}^K$ and scaling factors $\{\boldsymbol{\alpha}_k\}_{k=1}^K$ by applying downlink MAMSE receiver approach i.e., $\{\mathbf{U}_k \boldsymbol{\alpha}_k = \widetilde{\boldsymbol{\Gamma}}_k^{-1} \widehat{\mathbf{H}}_k^H \mathbf{G}_k \mathbf{P}_k\}_{k=1}^K$, where $\widetilde{\boldsymbol{\Gamma}}_k^{DL} = \widehat{\mathbf{H}}_k^H \mathbf{G} \mathbf{P} \mathbf{G}^H \widehat{\mathbf{H}}_k +$

$\sigma_{ek}^2 \text{tr}\{\mathbf{R}_{bk} \mathbf{G} \mathbf{P} \mathbf{G}^H\} \mathbf{R}_{mk} + \mathbf{R}_{nk}$. Note that in these expressions, $\{\boldsymbol{\alpha}_k\}_{k=1}^K$ are chosen such that each column of $\{\mathbf{U}_k\}_{k=1}^K$ has unity norm. Then, compute $\{\mathbf{B}_k, \mathbf{W}_k\}_{k=1}^K$ using (45).
> **Until** convergence.

**Convergence**: It can be shown that at each step of this algorithm, the sum AMSE of the system is non-increasing [3], [14]. Thus, the above iterative algorithm is guaranteed to converge. However, since all of our problems are non-convex, this iterative algorithm is not guaranteed to converge to the global optimum.

## X. EXTENSION OF **Algorithm II** TO OTHER PROBLEMS

In some cases, the power of each precoder entry needs to be constrained i.e,. $b_{ksn}^H b_{ksn} \leq \bar{\bar{p}}_{ksn}, \forall k, s, n$. This kind of power constraint has practical interest for coordinated BS systems where each BS is equipped with single antenna and the BS targets to allocate different powers to different symbols (see also [26], [27]). Mathematically, the latter problem can be formulated as

$$\mathcal{P}5: \min_{\{\mathbf{B}_k, \mathbf{W}_k\}_{k=1}^K} \bar{\xi}^{DL}, \text{ s.t } b_{ksn}^H b_{ksn} \leq \bar{\bar{p}}_{ksn}, \forall k, s, n. \quad (50)$$

It can be shown that this problem can be solved by **Algorithm II** with $\{\boldsymbol{\Delta}_{ks} = \text{diag}(\delta_{ks1}, \cdots, \delta_{ksN}), \forall s\}_{k=1}^K$. Consider the following problem variations

$$\mathcal{P}6: \min_{\{\mathbf{B}_k, \mathbf{W}_k\}_{k=1}^K} \sum_{s=1}^S p_s,$$
$$\text{s.t } \bar{\xi}^{DL} \leq \varepsilon_t, \text{ tr}(\sum_{k=1}^K \mathbf{B}_k \mathbf{B}_k^H) \leq P_{max} \quad (51)$$

$$\mathcal{P}7: \min_{\{\mathbf{B}_k, \mathbf{W}_k\}_{k=1}^K} \sum_{s=1}^S p_s,$$
$$\text{s.t } \bar{\xi}^{DL} \leq \varepsilon_t, [\sum_{k=1}^K \mathbf{B}_k \mathbf{B}_k^H]_{n,n} \leq \breve{p}_n, \forall n \quad (52)$$

$$\mathcal{P}8: \min_{\{\mathbf{B}_k, \mathbf{W}_k\}_{k=1}^K} \sum_{s=1}^S p_s,$$
$$\text{s.t } \bar{\xi}^{DL} \leq \varepsilon_t, \text{ tr}\{\mathbf{B}_k \mathbf{B}_k^H\} \leq \breve{p}_k, \forall k \quad (53)$$

$$\mathcal{P}9: \min_{\{\mathbf{B}_k, \mathbf{W}_k\}_{k=1}^K} \sum_{s=1}^S p_s,$$
$$\text{s.t } \bar{\xi}^{DL} \leq \varepsilon_t, \mathbf{b}_{ks}^H \mathbf{b}_{ks} \leq \bar{\bar{p}}_{ks}, \forall k, s \quad (54)$$

$$\mathcal{P}10: \min_{\{\mathbf{B}_k, \mathbf{W}_k\}_{k=1}^K} \sum_{s=1}^S p_s,$$
$$\text{s.t } \bar{\xi}^{DL} \leq \varepsilon_t, b_{ksn}^H b_{ksn} \leq \bar{\bar{p}}_{ksn}, \forall k, s, n \quad (55)$$

where $\varepsilon_t$ is the total sum AMSE target. The power allocation parts of $\mathcal{P}6 - \mathcal{P}10$ can be formulated as GPs. It is clearly seen that by modifying the power allocation step of $\mathcal{P}1$ to that of $\mathcal{P}6$, one can apply the solution approach of $\mathcal{P}1$ to solve $\mathcal{P}6$ (see also [4]). Next, we explain how the solution approach of $\mathcal{P}2$ can be extended to solve $\mathcal{P}7$. In the latter problem, each iteration should guarantee a non-increasing total

BS power. To ensure this non-increasing total BS power, (19) must be satisfied. This can be achieved just by modifying (24) to $\{\widetilde{\mathbf{b}}_n^H \widetilde{\mathbf{b}}_n \le \tilde{p}_n\}_{n=1}^N$. This leads to the following fixed point function.

$$\psi_n = \bar{\bar{f}}_n, \ \forall n \quad (56)$$

where

$$\bar{\bar{f}}_n = \frac{\tau}{\tilde{p}_n} \times$$
$$\frac{\psi_n [(\mathbf{A} + \boldsymbol{\Upsilon} + \boldsymbol{\Psi})^{-1}]_{(n,:)} \mathbf{A} ([(\mathbf{A} + \boldsymbol{\Upsilon} + \boldsymbol{\Psi})^{-1}]_{(n,:)})^H}{\sum_{i=1}^N \psi_i [(\mathbf{A} + \boldsymbol{\Upsilon} + \boldsymbol{\Psi})^{-1}]_{(i,:)} \mathbf{A} ([(\mathbf{A} + \boldsymbol{\Upsilon} + \boldsymbol{\Psi})^{-1}]_{(i,:)})^H}.$$

Therefore, one can apply the solution approach of $\mathcal{P}2$ to solve $\mathcal{P}7$ by modifying the power allocation step of $\mathcal{P}2$ to that of $\mathcal{P}7$ and computing $\{\psi_n\}_{n=1}^N$ with (56). The problems $\mathcal{P}8 - \mathcal{P}10$ can be solved like in $\mathcal{P}7$. The details are omitted for conciseness.

It is clearly seen that the analysis of this paper can be applied to solve other robust sum MSE-based constrained with groups of antennas, users or symbols power problems.

## XI. SIMULATION RESULTS

In this section, we present simulation results for $\mathcal{P}1 - \mathcal{P}4$. For all of our simulations, we take $K = 2$, $\{M_k = S_k = 2\}_{k=1}^K$ and $N = 4$. The entries of $\{\mathbf{R}_{bk}, \widetilde{\mathbf{R}}_{mk}\}_{k=1}^K$ are taken from a widely used exponential correlation model as $\{\mathbf{R}_{bk} = \rho_{bk}^{|i-j|}, \widetilde{\mathbf{R}}_{mk} = \rho_{mk}^{|i-j|}\}_{k=1}^K$, where $0 \le \rho_{bk}(\rho_{mk}) < 1$ and $1 \le i, j \le N(M_k), \forall k$. All of our simulation results are averaged over 100 randomly chosen channel realizations. It is assumed that $\sigma_{e1}^2 = 0.01, \sigma_{e2}^2 = 0.02, \rho_{b1} = 0.1, \rho_{b2} = 0.12, \rho_{m1} = 0.05, \rho_{m2} = 0.2, \mathbf{R}_{n1} = \sigma_1^2 \mathbf{I}_{M_1}$ and $\mathbf{R}_{n2} = \sigma_2^2 \mathbf{I}_{M_2}$. The Signal-to-Noise ratio (SNR) is defined as $P_{\text{sum}}/K\sigma_{av}^2$ and is controlled by varying $\sigma_{av}^2$, where $P_{\text{sum}}$ is the total BS power, $\sigma_2^2 = 2\sigma_1^2$ and $\sigma_{av}^2 = (\sigma_1^2 + \sigma_2^2)/2$. For the computation of SNR, we used the $P_{\text{sum}}$ obtained from the perfect CSI design of the proposed algorithm. We compare the performance of our iterative algorithm (**Algorithm II**) with that of [7] for the robust, non-robust and perfect CSI designs. The non-robust/naive design refers to the design in which the estimated channel is considered as perfect. Note that when $\{\sigma_{ek}^2 = 0\}_{k=1}^K$, the solution method of the latter paper turns to that of in [26].

### A. Simulation results for problem $\mathcal{P}1$

In this subsection, we compare the performance of our proposed algorithm with that of [7] when $P_{max} = P_{sum} = 10$. The comparison is based on the total sum AMSE and average symbol error rate (ASER)[4] of all users versus the SNR. Fig. 3.(a)-(b) show that **Algorithm II** and the algorithm in [7] achieve the same sum AMSE and ASER. Moreover, these figures show the superior performance of the robust design compared to that of the non-robust design. Next, we discuss the effects of antenna correlations factors $\{\rho_{bk}, \rho_{mk}\}_{k=1}^K$ on the system performance. For this purpose, we change the previous $\{\rho_{bk}, \rho_{mk}\}_{k=1}^K$ to $\rho_{b1} = 0.4, \rho_{b2} = 0.5, \rho_{m1} = 0.7$ and $\rho_{m2} = 0.8$ and plot the sum AMSEs of all designs in Fig. 3.(c). From Fig. 3.(a) and Fig. 3.(c), we can observe that the sum AMSE of non-robust, robust and perfect CSI designs increase as $\{\rho_{bk}, \rho_{mk}\}_{k=1}^K$ increase. This observation fits to that of [14] where $\mathcal{P}1$ is examined for $\{\mathbf{R}_{nk} = \sigma^2 \mathbf{I}\}_{k=1}^K$.

---

[4]For the sum AMSE design, ASER is also an appropriate metric for comparing the performance of the robust and non-robust designs [22]. QPSK modulation is utilized for each symbol.

### B. Simulation results for $\mathcal{P}2 - \mathcal{P}4$

In this subsection, we compare the performances of our proposed algorithm with that of [7] based on the total power utilized at the BS and the total achieved sum AMSE. We employ $\{\breve{p}_n = 2.5\}_{n=1}^N$ for $\mathcal{P}2$, $\{\breve{p}_k = 5\}_{k=1}^K$ for $\mathcal{P}3$, $\{\bar{\breve{p}}_{ks} = 2.5, \forall s\}_{k=1}^K$ for $\mathcal{P}4$ and $\{\rho_{bk}, \rho_{mk}\}_{k=1}^K$ as in the first paragraph of Section XI. As can be seen from Fig. 4.(a)-(c), for $\mathcal{P}2$, the proposed algorithm utilizes less total BS power than that of [7] for all designs. However, for $\mathcal{P}3$ and $\mathcal{P}4$, the proposed algorithm utilizes less total BS power than that of [7] for the robust designs only. Next, with the powers of Fig. 4.(a)-(c), we plot the total sum AMSEs of our algorithm and the algorithm in [7] as shown in Figs. 5.(a)-(c). These figures show that for the robust and non-robust designs, both of these algorithms achieve the same sum AMSE.

To examine the effects of antenna correlation factors, we use $\{\rho_{bk}, \rho_{mk}\}_{k=1}^K$ as in Section XI-A. Again Fig. 6.(a)-(c), show that our proposed algorithms utilize less total BS power than that of [7] in all designs and in the robust design only, for $\mathcal{P}2$ and $\mathcal{P}3 - \mathcal{P}4$, respectively. With the powers of Fig. 6.(a)-(c), we plot the sum AMSE of our algorithm and that of [7] as shown in Fig. 7.(a)-(c). These figures also show that both of these algorithms achieve the same sum AMSE.

In all figures of this section, the robust design outperforms the non-robust design and the improvement is larger for high SNR regions. This can be seen from the term $\boldsymbol{\Gamma}_k^{DL}$ of (2) where, at high SNR regions, $\sigma^2$ is negligible compared to $\sigma_{ek}^2 \text{tr}\{\mathbf{R}_{bk} \mathbf{B} \mathbf{B}^H\} \mathbf{R}_{mk}$ (the term due to channel estimation error). Thus, in the high SNR regions, since the non-robust design does not take into account the effect of $\sigma_{ek}^2 \text{tr}\{\mathbf{R}_{bk} \mathbf{B} \mathbf{B}^H\} \mathbf{R}_{mk}$ which is the dominant term, the performance of this design degrades. This implies that as the SNR increases, the performance gap between the robust and non-robust design increases. Furthermore, when $\{\rho_{bk}, \rho_{mk}\}_{k=1}^K$ increases, the system performance degrades. This is because as $\{\rho_{bk}, \rho_{mk}\}_{k=1}^K$ increases, the number of symbols with low channel gain increases (this can be easily seen from the eigenvalue decomposition of $\mathbf{R}_{bk}$ ($\mathbf{R}_{mk}$)). Consequently, for a given SNR value, the total sum AMSE also increases [14].

## XII. CONCLUSIONS

This paper considers sum MSE-based linear transceiver design problems for downlink multiuser MIMO systems where imperfect CSI is assumed to be available at the BS and MSs. These problems are examined for the generalized scenario where the constraint functions are per BS antenna, total BS, user or symbol power, and the noise vector of each MS is a ZMCSCG random variable with arbitrary covariance matrix. Each of these problems is solved as follows. First, we establish novel sum AMSE duality. Second, we formulate the power allocation part of the problem in the downlink channel as a GP. Third, using the duality result and the solution of GP, we utilize alternating optimization technique to solve the original downlink problem. We have established downlink-uplink duality to solve per BS antenna (groups of BS antenna) and total BS power constrained robust sum MSE-based problems. And, we have established downlink-interference duality to solve per user (groups of users) and symbol (groups of symbols) power constrained robust sum MSE-based problems. For the total BS power constrained robust sum MSE-based problems, the current duality is established by modifying the constraint function of the dual uplink channel problem. On the other hand, for the robust sum MSE minimization constrained with each BS antenna and each user (symbol) power problems, our duality are established by formulating the noise covariance matrices of the uplink and interference channels as fixed point functions, respectively. We have shown that our sum AMSE duality are able to solve any sum MSE-based robust design problem. Computer simulations verify the robustness of the proposed design compared to the non-robust/naive design.

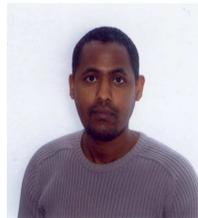

**Tadilo Endeshaw Bogale** (S'09) was born in Gondar, Ethiopia. He received his B.Sc and M.Sc degree in Electrical Engineering from Jimma University, Jimma, Ethiopia and Karlstad University, Karlstad, Sweden in 2004 and 2008, respectively. From 2004 to 2007, he was with the Mobile Project Department, Ethiopian Telecommunications Corporation (ETC), Addis Ababa, Ethiopia.

Since 2009 he has been working towards his PhD degree and as an Assistant Researcher with the ICTEAM institute, University Catholique de Louvain (UCL), Louvain-la-Neuve, Belgium. His research interests include robust (non-robust) transceiver design for multiuser MIMO systems, centralized and distributed algorithms, and convex optimization techniques for multiuser systems.

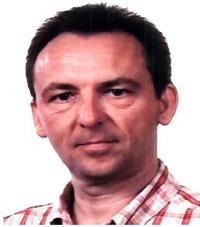
...

**Luc Vandendorpe** (M'93-SM'99-F'06) was born in Mouscron, Belgium, in 1962. He received the Electrical Engineering degree (summa cum laude) and the Ph.D. degree from the Universit Catholique de Louvain (UCL), Louvain-la-Neuve, Belgium, in 1985 and 1991, respectively. Since 1985, he has been with the Communications and Remote Sensing Laboratory of UCL, where he first worked in the field of bit rate reduction techniques for video coding. In 1992, he was a Visiting Scientist and Research Fellow at the Telecommunications and Traffic Control Systems Group of the Delft Technical University, The Netherlands, where he worked on spread spectrum techniques for personal communications systems. From October 1992 to August 1997, he was Senior Research Associate of the Belgian NSF at UCL, and invited Assistant Professor. He is currently a Professor and head of the Institute for Information and Communication Technologies, Electronics and Applied Mathematics.

His current interest is in digital communication systems and more precisely resource allocation for OFDM(A)-based multicell systems, MIMO and distributed MIMO, sensor networks, turbo-based communications systems, physical layer security and UWB based positioning.

Dr. Vandendorpe was corecipient of the 1990 Biennal Alcatel-Bell Award from the Belgian NSF for a contribution in the field of image coding. In 2000, he was corecipient (with J. Louveaux and F. Deryck) of the Biennal Siemens Award from the Belgian NSF for a contribution about filter-bank-based multicarrier transmission. In 2004, he was co-winner (with J. Czyz) of the Face Authentication Competition, FAC 2004. He is or has been TPC member for numerous IEEE conferences (VTC, Globecom, SPAWC, ICC, PIMRC, WCNC) and for the Turbo Symposium. He was Co-Technical Chair (with P. Duhamel) for the IEEE ICASSP 2006. He was an Editor for Synchronization and Equalization of the IEEE TRANSACTIONS ON COMMUNICATIONS between 2000 and 2002, Associate Editor of the IEEE TRANSACTIONS ON WIRELESS COMMUNICATIONS between 2003 and 2005, and Associate Editor of the IEEE TRANSACTIONS ON SIGNAL PROCESSING between 2004 and 2006. He was Chair of the IEEE Benelux joint chapter on Communications and Vehicular Technology between 1999 and 2003. He was an elected member of the Signal Processing for Communications committee between 2000 and 2005, and an elected member of the Sensor Array and Multichannel Signal Processing committee of the Signal Processing Society between 2006 and 2008. Currently, he is an elected member of the Signal Processing for Communications committee. He is the Editor-in-Chief for the EURASIP Journal on Wireless Communications and Networking. L. Vandendorpe is a Fellow of the IEEE.

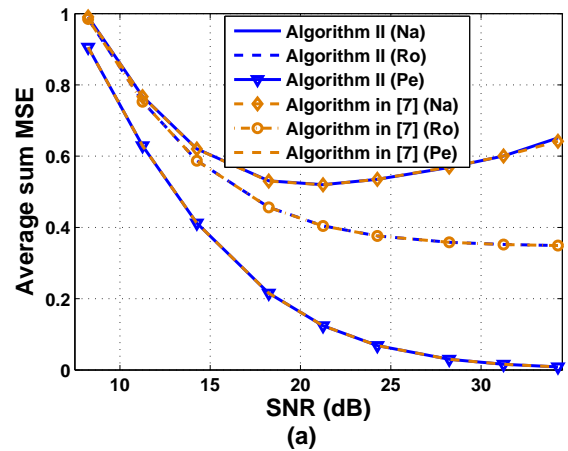

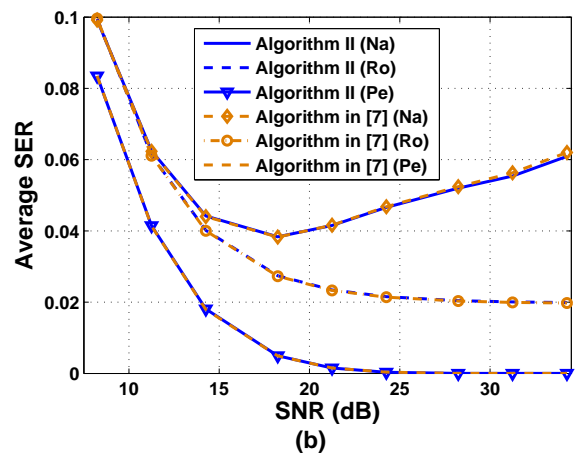

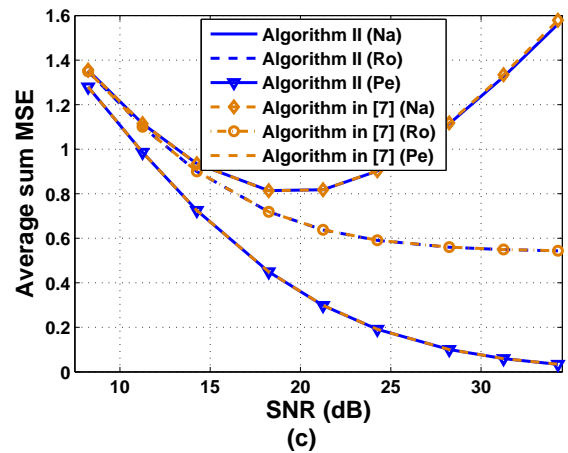

Fig. 3. Comparison of the proposed iterative algorithm (**Algorithm II**) and the algorithm in [7]. (a) In terms of total sum AMSE when $\rho_{b1} = 0.1, \rho_{b2} = 0.12, \rho_{m1} = 0.05$ and $\rho_{m2} = 0.2$. (b) In terms of ASER when $\rho_{b1} = 0.1, \rho_{b2} = 0.12, \rho_{m1} = 0.05$ and $\rho_{m2} = 0.2$. (c) In terms of total sum AMSE when $\rho_{b1} = 0.4, \rho_{b2} = 0.5, \rho_{m1} = 0.7$ and $\rho_{m2} = 0.8$. The non-robust/naive, robust and perfect CSI designs are denoted by (Na), (Ro), and (Pe), respectively.

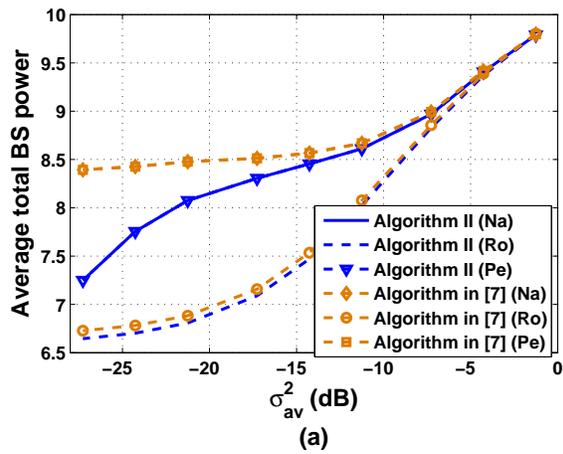
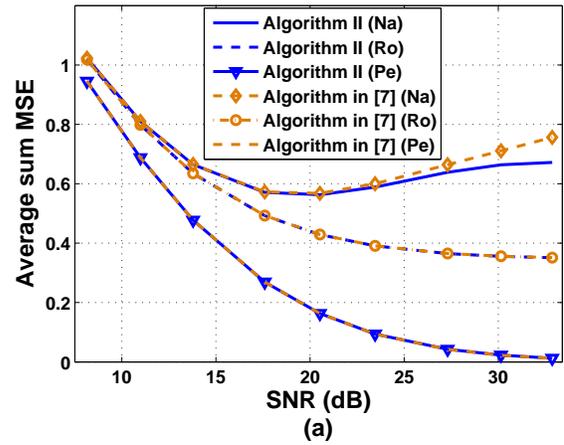
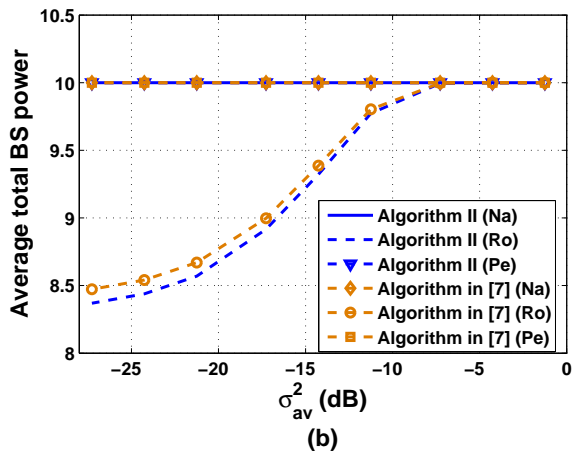
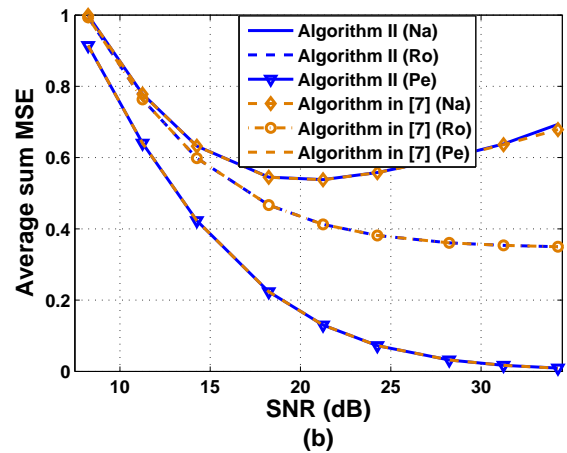
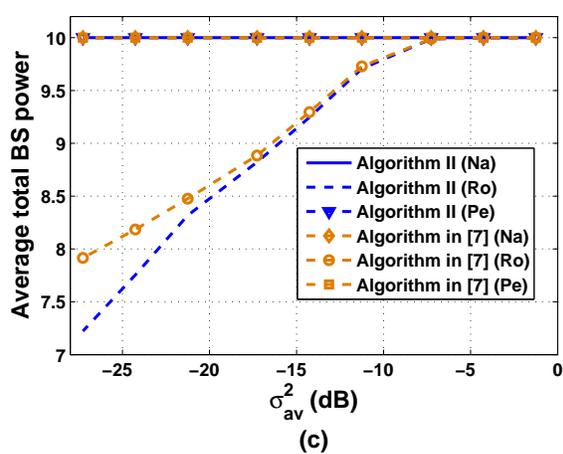
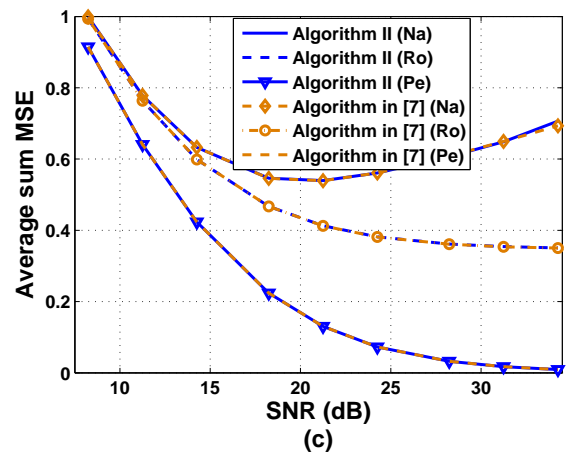

Fig. 4. Comparison of the proposed algorithm (**Algorithm II**) and the algorithm of [7] in terms of total BS power when $\rho_{b1} = 0.1, \rho_{b2} = 0.12, \rho_{m1} = 0.05$ and $\rho_{m2} = 0.2$. (a) for $\mathcal{P}2$. (b) for $\mathcal{P}3$. (c) for $\mathcal{P}4$.

Fig. 5. Comparison of the proposed algorithm (**Algorithm II**) and the algorithm of [7] in terms of total sum AMSE when $\rho_{b1} = 0.1, \rho_{b2} = 0.12, \rho_{m1} = 0.05$ and $\rho_{m2} = 0.2$. (a) for $\mathcal{P}2$. (b) for $\mathcal{P}3$. (c) for $\mathcal{P}4$.

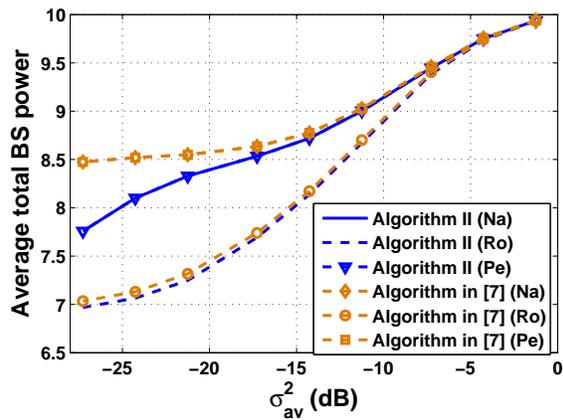
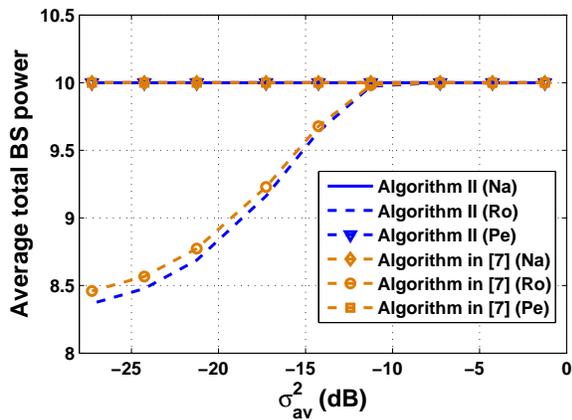
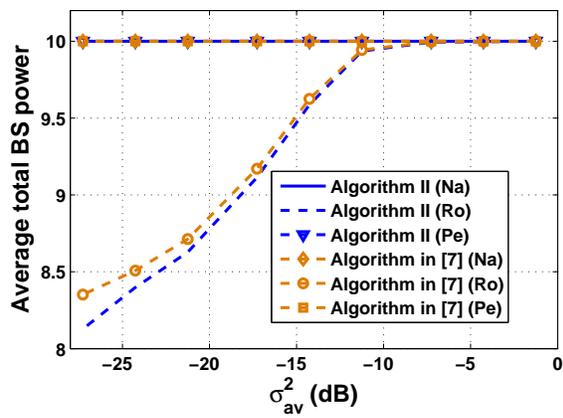

Fig. 6. Comparison of the proposed algorithm (**Algorithm II**) and the algorithm of [7] in terms of total BS power when $\rho_{b1} = 0.4, \rho_{b2} = 0.5, \rho_{m1} = 0.7$ and $\rho_{m2} = 0.8$. (a) for $\mathcal{P}2$. (b) for $\mathcal{P}3$. (c) for $\mathcal{P}4$.

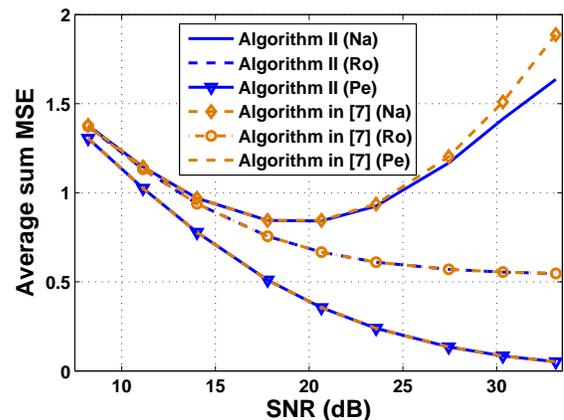
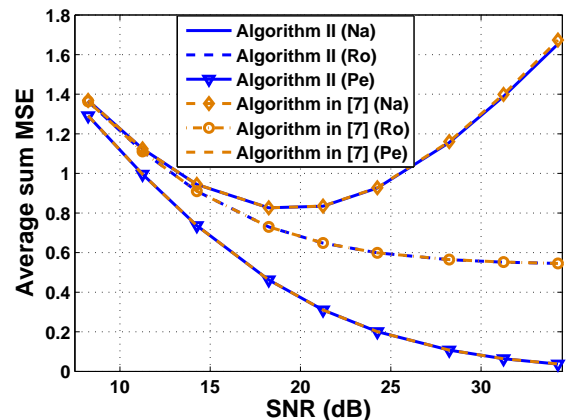
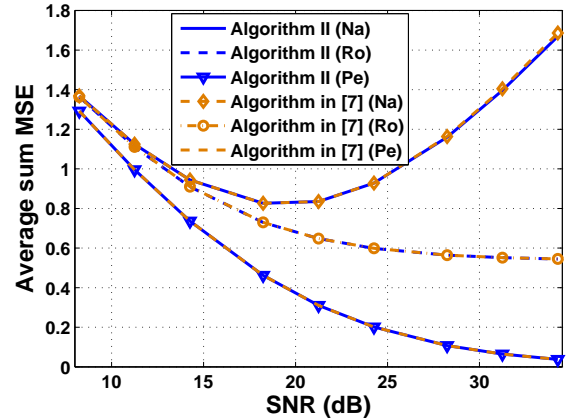

Fig. 7. Comparison of the proposed algorithm (**Algorithm II**) and the algorithm of [7] in terms of total sum AMSE when $\rho_{b1} = 0.4, \rho_{b2} = 0.5, \rho_{m1} = 0.7$ and $\rho_{m2} = 0.8$. (a) for $\mathcal{P}2$. (b) for $\mathcal{P}3$. (c) for $\mathcal{P}4$.